\documentclass[12pt]{iopart}

\usepackage[dvips]{epsfig}
\usepackage{amssymb,amsbsy,amsthm,amstext,amsfonts}
\usepackage{iopams}

\newtheorem{theorem}{Theorem}

\newtheorem{proposition}{Proposition}
\newtheorem{lemma}{Lemma}
\newtheorem{corollary}{Corollary}

\newcommand{\be}{\begin{equation}}
\newcommand{\ee}{\end{equation}}
\newcommand{\bea}{\begin{eqnarray}}
\newcommand{\eea}{\end{eqnarray}}
\newcommand{\beas}{\begin{eqnarray*}}
\newcommand{\eeas}{\end{eqnarray*}}
\newcommand{\nn}{\nonumber}
\newcommand{\bbm}{\begin{bmatrix}}
\newcommand{\ebm}{\end{bmatrix}}
\newcommand{\matl}{\left[ \begin{array}}
\newcommand{\matr}{\end{array} \right]}
\newcommand{\la}{\label}
\newcommand{\Lam}{\Lambda}
\newcommand{\La}{\mathcal{L}}

\newcommand{\bbr}{\mathbb{R}}
\newcommand{\mfrak}{\mathfrak}

\newcommand{\bS}{\mathbb{S}}
\newcommand{\wh}{\widehat}
\newcommand{\wt}{\widetilde}

\newcommand{\T}{^{\mbox{\small T}}}
\newcommand{\mdim}{\mbox{Dim }}
\newcommand{\mcal}{\mathcal}
\newcommand{\lin}{\langle\langle}
\newcommand{\rin}{\rangle\rangle}

\begin{document}

%------------------------------------------
% double space

%\setlength{\baselineskip}{15pt}

%-----------------------------------------------------------------

\title {A variational problem on Stiefel manifolds} 

\author{Anthony M. Bloch$^1$, Peter E. Crouch$^2$ and Amit K. Sanyal$^3$}
\address{$^1$ Department of Mathematics, 
University of Michigan, Ann Arbor, MI 48109}
\ead{abloch@math.lsa.umich.edu} 
\address{$^2$ Department of Electrical Engineering}
\address{$^3$ Department of Mechanical and Aerospace Engineering, 
Arizona State University, Tempe, AZ 85281}
\eads{\mailto{peter.crouch@asu.edu}, \mailto{sanyal@asu.edu}} 

\begin{abstract}
In their paper on discrete analogues of some classical systems such as the rigid 
body and the geodesic flow on an ellipsoid, Moser and Veselov introduced their 
analysis in the general context of flows on Stiefel manifolds. We consider here a 
general class of continuous time, quadratic cost, optimal control problems on Stiefel 
manifolds, which in the extreme dimensions again yield these classical physical %64
geodesic flows. We have already shown that this optimal control setting gives a
new symmetric representation of the rigid body flow and in this paper we extend this 
representation to the geodesic flow on the ellipsoid and the more general Stiefel 
manifold case. The metric we choose on the Stiefel manifolds is the same as that  
used in the symmetric representation of the rigid body flow, and that used by Moser 
and Veselov. In the extreme cases of the ellipsoid and the rigid body, the geodesic 
flows are known to be integrable. We obtain the extremal flows using both variational 
and optimal control approaches, and elucidate the structure of the flows on general 
Stiefel manifolds. %179 words
\end{abstract}

\ams{49Q99, 34K35, 65K10}
\submitto{Nonlinearity}

\maketitle

\section{Introduction}
  
This paper presents a variational problem on the Stiefel manifold of orthogonal $n$ 
frames in $N$ dimensional real Euclidean space, and its corresponding optimal control 
counterpart. Solutions to the variational and optimal control problems are obtained, 
and some of their geometric and analytic properties studied. This is an extension of 
earlier work by Bloch \etal (2002) on the symmetric representation of the rigid body 
equations (which correspond to the extreme case of $n=N$) on the Cartesian product 
$SO(N)\times SO(N)$. We characterize the space of solutions of the optimal control 
problem and the nature of the geodesic flows on the Stiefel manifold. 
The discrete version of this problem has been analyzed in the seminal work by Moser 
and Veselov (1991). Theorem 4 of Moser and Veselov (1991) gives a set of isospectral 
deformations for the discrete geodesic flow, which can be viewed as a discrete 
analogue of the parameter-dependent Lax representation. %These discrete isospectral 
%deformations give enough integrals for analytic integrability of the discrete dynamics. 
Recent results by Bolsinov and Jovanovic (2004) demonstrate that bi-invariant 
geodesic flows on Stiefel manifolds are integrable for a $SO(N)$-invariant metric. 
However, integrability has not yet been demonstrated for geodesic flows on general 
Stiefel manifolds with left-invariant metrics. \\
%We have not been able to determine the integrability of the corresponding continuous 
%geodesic flow. \\

We present a generalization of the Lax pair form for the equations
of motion and show how this reduces to the classical Lax pair form
in the case of the rigid body equations and the geodesic flow on the
ellipsoid. The integrability of the rigid body equations by a Lax pair formulation 
with parameter had been shown by Manakov (1976); Mischenko and Fomenko (1978) 
showed that a similar formalism exists for any semisimple Lie group. 
Further references on parameter-dependent Lax pair formulations for integrable 
systems are given in Fedorov (1995). The geodesic flow on an ellipsoid, which 
corresponds to the $n=1$ case, and its integrability have been treated by Moser and 
Veselov (1991), Knorrer (1980), and others. The paper by Moser and Veselov (1991) 
on the discrete variational version of this problem also gave a discrete Lax pair 
formulation with parameter, thereby demonstrating that the discrete geodesic 
equations on the Stiefel manifold are indeed integrable. \\

We describe here how to obtain geodesic flows using the maximum principle of optimal 
control theory (see Bloch \etal (2003), Gelfand and Fomin (2000), and Kirk (2004)). 
We show how to relate this optimal control formulation to the form naturally derived 
from variational calculus. We also relate these extremal flows to the Hamiltonian 
flow using the natural symplectic structure on the cotangent bundle of the Stiefel 
manifold. The extremal flows obtained here for the Stiefel manifolds are similar to 
the Hamiltonian flows on the ``extended Stiefel varieties" as described in Federov 
(2005). Finally, we demonstrate that the natural symplectic manifold carrying the 
extremals of the optimal control problem is symplectomorphic to the cotangent bundle 
of the Stiefel manifold. In the next section, we pose our problem on a general 
Stiefel manifold, and give the extremal flows obtained in the limiting cases of the 
sphere/ellipsoid ($n=1$), and the $N$ dimensional rigid body ($n=N$). Section 3 
presents the extremal solution to the variational problem for the general case, and 
connects this solution to the extremal flow on the cotangent bundle, given by 
equation (\ref{VPsol}). Section 4 gives the extremal solution of the optimal control 
problem restricted to a symplectic submanifold of $R^{nN}\times R^{nN}$ that has the 
dimension of the cotangent bundle. This section also gives the correspondence between 
this extremal solution and the extremal solution of the variational problem. Section 
5 presents the structure of the tangent and cotangent bundles of the Stiefel manifold, 
and establishes a symplectomorphism between the manifold carrying the extremal 
solutions of the optimal control problem, and the cotangent bundle. Section 6 
presents a few applications and some future research issues of interest regarding 
geodesic flows on Stiefel manifolds, while Section 7 presents some concluding remarks.

\section{Background and Limiting cases}

We introduce the variational and optimal control problems on a Stiefel manifold in 
this section, based on minimizing the time integral of the kinetic energy. The 
metric on the manifold is given by the kinetic energy expression. We also give 
the extremal flows obtained in the limiting cases of the sphere/ellipsoid 
($n=1$), and the $N$ dimensional rigid body ($n=N$). The extremal flows in these 
cases are well-known and integrable, and have been given in several earlier works 
such as Knorrer (1982), Moser (1980), and Bloch \etal (2002).

\subsection{Variational and Optimal Control Problems on a Stiefel Manifold}

The Stiefel manifold $V(n,N)\subset\bbr^{nN}$ consists of orthogonal $n$ frames in 
$N$ dimensional real Euclidean space,
\[ V{(n,N)}=\{Q\in\bbr^{nN};\quad QQ\T=I_n\}. \]
Introduce the pairing in $\bbr^{rs}$ given by 
\be \langle A,B\rangle =\Tr (A\T B), \la{pair} \ee
where $\Tr (\cdot)$ denotes trace of a matrix and the left invariant metric on 
$\bbr^{nN}$ given by 
\be \lin W_1,W_2\rin = \langle W_1\Lam, W_2\rangle= \langle W_1,W_2\Lam\rangle, 
\la{metric} \ee
where $\Lam$ is a positive definite $N\times N$ diagonal matrix. The pairing 
(\ref{pair}) was used in Ratiu (1980) as a positive definite bilinear form on 
$\mfrak{so}(n)$. We are interested in the variational problem given by: 
\be
\min_{Q(\cdot)}\int^T_0\frac12\lin\dot Q,\dot Q\rin dt
\label{VPnN} \ee
%\hfill VP(n,N)\newline
\noindent
subject to: $QQ\T=I_n$, $Q\in\bbr^{nN}$, $1\le n\le N$, $Q(0)=Q_0$, $Q(T)=Q_T$, $I_n$ 
denotes the $n\times n$ identity matrix. This is a variational problem defined on the 
Stiefel manifold $V(n,N)$. The dimension of this manifold is given by 
\[ \mbox{Dim }V(n,N)=nN-\frac{n(n+1)}{2}=n(N-n)+\frac{n(n-1)}{2}. \] 
The corresponding optimal control problem is given by:
\be
\min_{U(\cdot)}\int^T_0\frac12\lin QU,QU\rin dt
\la{OCPnN} \ee
%\hfill OCP(n,N)\newline
\noindent
subject to: $\dot Q=QU$; $QQ\T=I_n$, $Q(0)=Q_0$, $Q(T)=Q_T$ where $U\in \mfrak{so}(N)$. 
Note that the quantity to be minimized is invariant with respect to the left action of 
$SO(n)$ on $V(n,N)$ since the metric (\ref{metric}) is left invariant. 

\subsection{The Rigid Body equations}

For the special case when $n=N$, $V(N,N)\equiv SO(N)$ and the extremal trajectories 
of the optimal control problem (\ref{OCPnN}) give the $N$-dimensional rigid body 
equations. The usual system of rigid body equations on $T^\star SO(N)$ are 
\be \dot Q= QU,\;\ \ \dot M= [M,U], \la{TsSON} \ee
where $M=J(U)\triangleq \Lam U+U\Lam$ is the body momentum and $[\cdot,\cdot]$ denotes 
the matrix commutator. Equations (\ref{TsSON}) 
can also be obtained directly from the variational problem (\ref{VPnN}) in the case 
$n=N$. As shown by Manakov (1976) and Ratiu (1980), the Euler equations 
are integrable with a parameter-dependent Lax pair representation 
\[ \frac{d}{d t} (M+\lambda \Lam^2)= [M+\lambda\Lam^2,U+\lambda\Lam]. \]
In the optimal control approach, the costate $P$ is a vector of Lagrange multipliers 
used to enforce the equality constraint $\dot Q=QU$ (see Gelfand and Fomin, 2000). 
The extremal trajectories for the optimal control problem in this case are given by 
\be \dot Q= QU,\;\ \ \dot P= PU, \la{symRB} \ee
where $U\in\mfrak{so}(N)$ and $[\cdot,\cdot]$ denotes the matrix commutator.
In the symmetric representation of the rigid body 
equations given in Bloch \etal (2002), the states $Q$ and the costates $P$ are both 
orthogonal $N\times N$ matrices, and the extremal trajectories are on $SO(N)\times 
SO(N)$. For the symmetric representation of the rigid body equations (\ref{symRB}), 
$U$ is regarded as a function of $Q$ and $P$ via the equations 
\[ U\triangleq J^{-1}(M)\in\mfrak{so}(N),\;\ \mbox{ and }\;\ M\triangleq Q\T P-P\T
Q. \] 
The inverse of the mapping $\Phi: (Q,P) \mapsto (Q,M)=(Q,Q\T P-P\T Q)$ is given by 
\be P=Q\big(e^{\sinh^{-1}\frac M2}\big), \la{PMrel} \ee
where $\mbox{sinh: }\mfrak{so}(N)\rightarrow \mfrak{so}(N)$ and its inverse are 
defined in Bloch \etal (2002). The spatial momentum given by 
\[ m=QMQ\T = PQ\T-QP\T, \]
is conserved along the flow given by (\ref{symRB}). Equivalence between (\ref{symRB}) 
and (\ref{TsSON}) was established by Bloch \etal (2002) on the sets $S$ and $S_M$ 
where 
\beas S&=&\{(Q,P)\in SO(N)\times SO(N);\quad \|M\|_{OP}<2\}, \\
S_M&=&\{(Q,M)\in SO(N)\times \mfrak{so}(N);\quad\|M\|_{OP}<2\}, \eeas 
and $\|A\|_{OP}= \mbox{sup}\{\|Ax\|\ :\ \|x\|=1\}$ is the operator norm. $\Phi^{-1}$ 
is then well defined on $S_M$. \\

Note that for the extremal flow to the optimal control problem given by (\ref{symRB}), 
one need not choose the costates to be orthogonal. One can choose a costate vector 
$P_0$ such that $Q\T P_0$ is skew symmetric, for example. In that case, the extremal 
flow is given by 
\be \dot Q= Q J^{-1}(Q\T P_0), \quad \dot P_0= P_0 J^{-1}(Q\T P_0), \la{macsco} \ee 
and these equations restrict to the invariant submanifolds defined by $Q\T P_0\in 
\mfrak{so}(N)$. These are the McLachlan-Scovel equations (see McLachlan and Scovel, 
1995). Comparing these equations with (\ref{symRB}), we see that $P_0=QM= P-QP\T Q$. 
In Section \ref{sec:ocpro}, we obtain a generalization of the McLachlan-Scovel 
equations to the Stiefel manifolds $V(n,N)$, where $1<n<N$. 

\subsection{Geodesic flow on the ellipsoid}

For the other extreme case, when $n=1$, we obtain the equations for the geodesic flow 
on the sphere $V(1,N)\equiv \bS^{N-1}$ with $Q=q\T$, $q\T q=1$. This can be also be 
regarded as the geodesic flow on the ellipsoid 
\[ \bar q\T\Lam^{-1}\bar q=1, \]
where $q=\Lam^{-1/2}\bar q$. The costate variable $P=p\T$ is used to enforce the 
constraint $\dot q= -Uq$ for the optimal control problem (\ref{OCPnN}) when $n=1$. 
The extremal solutions to this optimal control problem are
\be \dot q= -Uq,\;\;\ \dot p= -Up+Aq, \la{SNp1} \ee
where $A=qq\T U\Lam U- U\Lam U qq\T$. These extremal solutions have the same form as
the extremal solutions for the general case (when $1<n<N$) given in Section 3. 
The body momentum is obtained as   
\be M=qp\T - pq\T, \la{bMomel} \ee
in terms of the solution $(q,p)$ to the optimal control problem. Equations (\ref{SNp1} 
can than be expressed in terms of the body momentum as 
\be \dot q= -Uq,\;\; \dot M= [M,U]-A. \la{bMomeq} \ee
The body momentum can also be expressed in terms of the solution to the variational 
problem $(q,\dot q)$, which we write in the form $(q,s)$, where 
$s=\Lam\dot q=-\Lam Uq$, as 
\[ M=qs\T -sq\T= qq\T U\Lam +\Lam Uqq\T. \] 
Pre- and post-multiplying the above expression with $\Lam^{-1}$, we obtain
\be \Lam^{-1} M\Lam^{-1}= \Lam^{-1} qq\T U+ Uqq\T\Lam^{-1}. \la{LML} \ee
Since $U$ is skew-symmetric, $q\T Uq=0$, and post-multiplying both sides of equation 
(\ref{LML}) with the vector $q$, we get
\be \Lam^{-1} M\Lam^{-1} q = Uq(q\T\Lam^{-1} q) \;\
\Leftrightarrow Uq = \frac{\Lam^{-1} M\Lam^{-1}}{q\T\Lam^{-1} q}q. \la{Uqe} \ee
Note that equation (\ref{Uqe}) specifies $U$ upto the equivalence class 
\[ [U],\; V\in [U] \Leftrightarrow Uq=Vq. \] 
Equation (\ref{Uqe}) gives our choice for $U$ as 
\[ U= \frac{\Lam^{-1} M\Lam^{-1}}{q\T\Lam^{-1} q}. \]
From this expression for $U$, we obtain 
\[ \fl U\Lam U= \frac1{(q\T\Lam^{-1} q)^2} \Lam^{-1}\Big(p\T\Lam^{-1} q(qp\T+ pq\T)-
(q\T\Lam^{-1} q)pp\T-(p\T\Lam^{-1} p)qq\T\Big)\Lam^{-1}, \]
and \[ qq\T U\Lam U= \frac1{(q\T\Lam^{-1} q)^2} qq\T\Lam^{-1}\Delta, \]
where \[ \Delta= (p\T\Lam^{-1} q)^2- (p\T\Lam^{-1} p)(q\T\Lam^{-1} q). \]
This gives us $A$ in equation (\ref{SNp1}) as 
\be A=[qq\T, U\Lam U]= \frac1{(q\T\Lam^{-1} q)^2}\big(qq\T\Lam^{-1}-\Lam^{-1} qq\T\big)
\Delta. \la{Aeqn} \ee 
Equations (\ref{SNp1}), (\ref{Uqe}) and (\ref{Aeqn}) defines extremal trajectories 
of the optimal control problem on $V(1,N)=\bS^{N-1}$. \\

The Lagrangian (variational) formulation for this problem gives us the equations for 
the geodesic flow on the sphere.
%\be \dot q= -Uq,\;\;\ \dot s= -Us+Aq. \la{SNp1v} \ee 
%Applying Proposition 2, 
To obtain these equations, we take reduced variations (see Marsden and Ratiu, 1999) on 
$V(1,N)=\bS^{N-1}$. The equation of motion can be written as 
\be \Lam \ddot q= bq, \la{ddq} \ee
where $b$ is a real scalar in this case. We determine $b$ from the constraint $q\T q=1$. 
Differentiating this constraint with respect to time twice, we get 
\[ \ddot q\T q+ \dot q\T\dot q= 0. \]
Substituting for $\ddot q$ from equation (\ref{ddq}) into the above expression, we get 
\[ b= -\frac{\dot q\T\dot q}{q\T\Lam^{-1} q}. \]
Thus, we get the Lagrangian (variational) equations for the geodesic flow on the 
sphere ($\bS^{N-1}$) as 
\be \ddot q= -\frac{\dot q\T\dot q}{q\T\Lam^{-1} q}\Lam^{-1} q. 
\la{varSNm1} \ee 
Integrability of these extremal flows were proven by Jacobi with relation to the 
Neumann problem of motion on the sphere with a quadratic potential, as shown by 
Knorrer (1982). Contemporary version of integrability of the geodesic flow on an 
ellipsoid was demonstrated by Moser (1980) using the Theorem of Chasles and 
geometry of quadrics. We now obtain a Manakov Lax pair formulation for this system. 
In this case, we have 
\[ \Lam U-U\Lam= (M\Lam^{-1}-\Lam^{-1} M)/q\T\Lam^{-1} q. \]
Using the body momentum equation in (\ref{bMomeq}), we now get  
\be \frac{d}{dt}M(\lambda)= [M(\lambda),U(\lambda)], \la{kock} \ee
where
\bea M(\lambda) &=& M+\lambda\Lam+qq\T\lambda^{-1}\Delta/q\T\Lam^{-1}q, \\
U(\lambda) &=& U-\lambda\Lam^{-1}/q\T\Lam^{-1}q. \la{Ulam} \eea 
From (\ref{kock}), we see that the coefficients of $\lambda^j$ in the expansion of 
$\frac{1}{2k}\mbox{tr}(M(\lambda))^k$ are conserved along the geodesic flow. 
%, and give enough integrals in involution to ensure integrability. 

\section{The Variational Problem}\la{sec:varp}

Now we look at the solution to the general problem (\ref{VPnN}) posed on the Stiefel 
manifold $V(n,N)$ where $1\le n\le N$. The Lagrangian for this variational problem is 
\be \La (Q,\dot Q)= \frac12\lin\dot Q,\dot Q\rin = \frac12 \Tr(\Lam\dot Q\T \dot Q), 
\la{Lag} \ee 
where $Q\in V(n,N)$, i.e., $Q$ satisfies $QQ\T=I_n$.

\subsection{Solution to the Variational Problem}

The solution to the variational problem (\ref{VPnN}) posed on $V(n,N)$ is given by 
the following result.
\begin{theorem}
The variational problem (\ref{VPnN}) on the Stiefel manifold $V(n,N)$ has the solution: 
\bea \lefteqn{\dot Q= QU, \;\ U\in \mfrak{so}(N),} \la{Qdot} \\ 
\lefteqn{Q\T Q\dot U \Lam+ \Lam\dot U Q\T Q= \Lam U^2 Q\T Q-Q\T Q U^2\Lam.} \la{solVP} 
\eea
\la{vaSol}
\end{theorem}
\noindent{\em Proof:} We take reduced variations on the Stiefel manifold 
given by $QQ\T=I_n$; this ensures that the variation vector field is always locally 
tangent to the manifold. The reduced variations are given by 
\be \delta Q= Q\Sigma,\;\ \delta\dot Q= \dot Q \Sigma+ Q\dot\Sigma, \la{varns} \ee
where $Q\Sigma\in T_Q V(n,N)$ and $\Sigma$ is an $N\times N$ skew-symmetric matrix. 
These variations have fixed end points, i.e., $\Sigma(0)=\Sigma(T)=0$. The kinematic 
expression (\ref{Qdot}) is obtained easily from the constraint $QQ\T=I_n$. On taking 
the first variation of the integral quantity in (\ref{VPnN}) and setting it to zero, 
we obtain 
\[ \int^T_0 \{ \langle \delta\dot Q\Lam,\dot Q\rangle+ \langle\dot Q\Lam, 
\delta\dot Q\rangle \}dt=0. \]
On carrying out integration by parts using fixed end point variations, the above 
expression simplifies to 
\bea - \int^T_0 \{\langle Q\Sigma\Lam,\ddot Q\rangle+ \langle\ddot Q\Lam, Q
\Sigma\rangle \}dt &=& 0 \nn \\
\Leftrightarrow \int^T_0 \langle\ddot Q\Lam, Q\Sigma\rangle dt &=& 0 \nn \\
\Leftrightarrow \Tr (\Lam\ddot Q\T Q\Sigma) &=& 0 \nn \\
\Leftrightarrow \Lam\ddot Q\T Q- Q\T\ddot Q\Lam &=& 0, \la{vareq} \eea
i.e., $\Lam\ddot Q\T Q$ is symmetric. Taking the time derivative of $\dot Q$ in 
(\ref{Qdot}), we get
\[ \ddot Q= QU^2+Q\dot U. \]
Substituting this for $\ddot Q$ in equation (\ref{vareq}) gives the solution of the 
variational problem as given by equation (\ref{solVP}). \qed

Note that the flow on the Stiefel manifold induced by equation (\ref{solVP}) 
is invariant to variations within the equivalence class $[U]\subset \frak{so}(N)$, 
defined by 
\be V\in [U]\ \Leftrightarrow\ QV=QU. \la{eqvclass} \ee
Since $QQ\T= I_n$, the left-multiplication of equation (\ref{solVP}) by $Q$ gives
\be Q\dot U\Lam Q\T+ Q\Lam\dot U Q\T= Q[\Lam,U^2]Q\T. \la{vasol} \ee
Here $W=\dot U\Lam+ \Lam\dot U \in \frak{so}(N)$ is a skew-symmetric matrix that also 
satisfies $QWQ\T=Q[\Lam, U^2]Q\T$. Equation (\ref{vasol}) is a statement of the 
conservation of the quantity 
\be m= QU\Lam Q\T+ Q\Lam U Q\T= SQ\T-QS\T, \la{consm} \ee
where $S=QU\Lam=\dot Q \Lam$, along the flow of the vector field given by equations 
(\ref{Qdot})-(\ref{solVP}). Note that 
\[ \dot m = QWQ\T- Q[\Lam, U^2]Q\T= 0 \] 
by equation (\ref{vasol}). The quantity $m$ is analogous to the spatial momentum of 
the rigid body, which is conserved along the flow of the rigid body equations. The 
equations (\ref{Qdot})-(\ref{solVP}) can also be written as 
\bea \lefteqn{\dot Q= QU= S\Lam^{-1}, \;\ U\in \mfrak{so}(N),} \la{soVP} \\ 
\lefteqn{Q\T \dot S- \dot S\T Q= 0.} \la{sVP} \eea 
%The above equations are a compact form for equations (\ref{Qdot})-(\ref{solVP}), 
%and express the extremal solution to this variational problem in the Lagrangian 
%formulation $(Q,\dot Q)$ with the pair $(Q,S)$. 

The above formulation of the variational problem, equations (\ref{soVP})-(\ref{sVP}) 
in particular, do not give an implicit or explicit equation for the rate $\dot S$. 
To obtain such an expression, we need the following result. 
\begin{proposition}
Define the linear map $L_Q:\ \bbr^{nN} \rightarrow \bbr^{NN}$, given by 
\be L_Q(X)= Q\T X- X\T Q, \;\ Q\in V(n,N). \la{LQmap} \ee
If $X\in \bbr^{nN}$ is in the kernel of this map, then $X$ is of the form %satisfies 
%\be Q\T X= X\T Q, \la{LQ} \ee and $X$ is of the form 
\be X= BQ,\;\ \ B=B\T. \la{kern} \ee 
\la{LQker} 
\end{proposition}
\vspace*{-3mm}
\noindent {\em Proof}: If $X$ is in the kernel of the $L_Q$ defined by (\ref{LQmap}), 
then $Q\T X=X\T Q$ and hence $Q\T X$ is symmetric. Let us denote the rows of $Q$ by 
the orthonormal vectors $q_1\T$, $q_2\T$, $\ldots$, $q_n\T$. We complete a 
right-handed orthonormal basis of $\bbr^N$ from this set of orthonormal vectors; the  
additions to this set are denoted by the vectors $q_{n+1}\T,\ldots,q_N\T$. The 
vectors of this orthonormal basis are arranged in columns to form the orthogonal 
matrix 
\[ C= \left[ q_1\;\ q_2\;\ \cdots\;\ q_n\;\ q_{n+1}\;\ \cdots\;\ q_N \right] \in 
SO(N). \]
Thus, on right-multiplication of $Q$ by $C$ we get: 
\be R=Q C= [I_n\;\ 0]. \la{QSI} \ee
We express $C\T$ in the block diagonal form 
\[ C\T= \matl{cc} C_1 & C_3 \\ C_2 & C_4 \matr, \] 
where the diagonal blocks $C_1$ and $C_4$ are $n\times n$ and $(N-n)\times (N-n)$ 
matrices respectively. Hence, from equation (\ref{QSI}), we get
\be Q= [I_n\;\ 0]C\T= [C_1\;\ C_3]. \la{S1S3} \ee
We apply the same transformation to $X\in \bbr^{nN}$ to get 
\be Y= XC= [Y_1\;\ Y_2], \la{CQSig} \ee
where $Y_1\in\bbr^{nn}$ and $Y_2\in\bbr^{n(N-n)}$. Since $Q\T X= CR\T YC\T$, the 
relation $Q\T X=X\T Q$ can be expressed as 
\be R\T Y= Y\T R\ \Leftrightarrow\ Y_1=Y_1\T=B,\;\ Y_2=0. \la{Yev} \ee
Then we have 
\beas X&=& [Y_1\;\ 0] \matl{cc} C_1 & C_3 \\ C_2 & C_4 \matr \\
&=& [BC_1\;\ BC_3] \\
&=& B[C_1\;\ C_3]= BQ, \eeas
using the expression in equation (\ref{S1S3}) for $Q$. This gives us the result in 
equation (\ref{kern}), where $B=Y_1=Y_1\T$ is a symmetric $n\times n$ matrix. \qed \\

Since from (\ref{sVP}), $\dot S$ is in the kernel of $L_Q$, using the above result 
we get 
\be \dot S= BQ, \;\ B=B\T. \la{Sd} \ee
The following lemma is necessary to obtain an expression for the flow of the 
quantity $S$ and leads to an implicit equation for $B$. 
\begin{lemma}
If $K$ is symmetric and positive definite and $B$ is symmetric, the map $J: 
\mbox{Sym}(n)\rightarrow \mbox{Sym}(n)$ given by $J: B \mapsto KB+BK$ has kernel 
zero, and is hence an isomorphism. 
\la{JBlem} 
\end{lemma}
\noindent{\em Proof}: Since $B$ is symmetric, there exists an orthogonal matrix $C$, 
such that $CBC\T= \Psi$, where $\Psi$ is a real diagonal matrix. Thus if $KB+BK=0$, 
then 
\[ C\T KC\Psi+ \Psi C\T KC=0, \]
where $C\T KC$ is a positive definite symmetric matrix. If $e$ is an eigenvector 
of $\widehat{K}= C\T KC$, then 
\[ \widehat{K} e= \lambda e,\quad \lambda>0, \]
and 
\[ \widehat{K}\Psi+ \Psi\widehat{K}=0 \Rightarrow \widehat{K}(\Psi e)+\lambda 
(\Psi e)= 0. \]
Hence $\Psi e$ is also an eigenvector of $\widehat{K}$ with eigenvalue $-\lambda 
<0$. But $\widehat{K}$ is positive definite and so all its eigenvalues are strictly 
positive. Thus, we have a contradiction, unless $\Psi=0$ and hence $B=0$. \qed \\ 
Here $\mbox{Sym}(n)$ denotes the space of $n\times n$ symmetric matrices. Now we
state the main result of this subsection, which expresses the geodesic flows on 
the Stiefel manifold in terms of the quantities $Q$ and $S=QU\Lam$; the pair  
$(Q,S)$ can be used to parametrize the tangent bundle since $S=\dot Q\Lam^{-1}$.  
\begin{proposition}
The geodesic flow on the Stiefel manifold given by the variational problem 
(\ref{VPnN}) is of the form 
\be \dot Q= S\Lam^{-1},\quad \dot S= BQ, \la{QSeq} \ee
where 
\be B=\wt{J_Q}^{-1}(E),\quad E= -2\dot Q\dot Q\T= -2S\Lam^{-2} S\T, \la{Bexp} \ee 
and 
\be \wt{J_Q}(B)\triangleq (Q\Lam^{-1}Q\T) B+ B(Q\Lam^{-1} Q\T);\;\ \wt{J_Q}: 
\mbox{Sym}(n) \rightarrow \mbox{Sym}(n). \la{JQB} \ee 
\la{varexs}
\end{proposition}
\noindent{\em Proof}: The proof of this result makes use of the simple observation 
that 
\[ Q\Lam^{-1}S\T= -QUQ\T \]
is skew-symmetric. We already know from Proposition \ref{LQker} and equation 
(\ref{sVP}) that $\dot S= BQ$ where $B$ is a $n\times n$ symmetric matrix. To
obtain an expression for $B$, we take a time derivative of 
\[ Q\Lam^{-1}S\T+ S\Lam^{-1}Q\T=0 \]
along the extremal trajectories. This gives us 
\beas Q\Lam^{-1}Q\T B+ BQ\Lam^{-1}Q\T&=& S\Lam^{-1}UQ\T-QU\Lam^{-1}S\T \\ 
\Rightarrow \wt{J_Q}(B)&=& -2\dot Q\dot Q\T= -2S\Lam^{-2}S\T, \eeas 
since $S=QU\Lam$. Since $Q\Lam^{-1}Q\T$ is positive definite, we know from Lemma 
\ref{JBlem} that $\wt{J_Q}: \mbox{Sym}(n) \rightarrow \mbox{Sym}(n)$ is an 
isomorphism, and hence the inverse $B=\wt{J_Q}^{-1}(2QU^2 Q\T)$ exists and is 
unique. This proves the proposition. \qed \\
Note that the quantity $S=QU\Lam$ depends on the equivalence 
class defined by (\ref{eqvclass}); so equations (\ref{QSeq})-(\ref{JQB}) in terms 
of $(Q,S)$ also have the advantage that they {\em uniquely} express the extremal 
flows, whereas $U$ in equations (\ref{Qdot})-(\ref{solVP}) is not unique. Since 
$S=\dot Q\Lam^{-1}$, the above proposition also expresses the extremal flow in 
terms of the tangent bundle pair $(Q,\dot Q)\in TV(n,N)$. \\ 

Now we give another method to obtain the above geodesic flow from the constrained 
Lagrangian 
\be \La_c (Q,\dot Q)= \La (Q,\dot Q)+ \frac12\langle B,QQ\T-I_n\rangle= \frac12\lin 
\dot Q,\dot Q\rin+ \frac12\langle B,QQ\T-I_n\rangle, \la{consLag} \ee 
where $B\in\mbox{Sym}(n)$ is a Lagrange multiplier matrix. We form the Hamiltonian 
for this problem by applying the Legendre transform to the Lagrangian (\ref{Lag}), 
which gives the same result as applying the Legendre transform to the constrained 
Lagrangian (\ref{consLag}). Just as $T V(n,N)$ can be viewed as a submanifold of 
$V(n,N)\times\bbr^{nN}\subset \bbr^{nN}\times\bbr^{nN}$, $T^\star V(n,N)$ can also 
via $\langle\cdot,\cdot\rangle$ be viewed as a submanifold of $V(n,N)\times\bbr^{nN}
\subset\bbr^{nN}\times\bbr^{nN}$, and endowed with the symplectic structure induced 
from the canonical symplectic structure of $\bbr^{nN}\times\bbr^{nN}$. 
\begin{proposition} \la{legend}
The Legendre transform $\mathbb{F}\mcal{L}\ :\ TV(n,N)\rightarrow T^\star V(n,N)$ 
is obtained from the fiber derivative of the Lagrangian $\La\ :\ TV(n,N)\rightarrow 
\bbr$ as 
\be \mathbb{F}\mcal{L} (Q,\dot Q)= S. \la{legt} \ee
\end{proposition}
\noindent {\em Proof}: Let $U,\ V\in\mfrak{so}(N)$, so that $QU,\ QV\in T_Q V(n,N)$. 
The fiber derivative of the Lagrangian (\ref{Lag}) is given by 
\beas \mathbb{F}\mcal{L}(Q,\dot Q) (\Gamma) &=& \frac{d}{ds}\big|_{s=0} 
\La (\dot Q+s\Gamma) \\ &=& \frac{d}{ds}\big|_{s=0} \frac12\lin \dot Q+s\Gamma, 
(\dot Q+s\Gamma)\rin \\ &=& \frac12\lin \Gamma,\dot Q\rin+ \frac12\lin \dot Q,
\Gamma\rin = \lin \dot Q,\Gamma\rin \\ 
\Leftrightarrow \mathbb{F}\mcal{L}(Q,\dot Q) (\Gamma) &=& \langle QU\Lam,\Gamma\rangle
=\langle S,\Gamma\rangle,  \eeas 
or by identifying $T^\star V(n,N)$ with $T V(n,N)$ through $\langle\cdot,\cdot\rangle$, 
$\mathbb{F}\mcal{L}(Q,\dot Q)= S$. Thus, the map given by (\ref{legt}) is the Legendre 
transform $\mathbb{F}\mcal{L}\ :\ TV(n,N)\rightarrow T^\star V(n,N)$. \qed \\
Proposition \ref{legend} can be used to 
obtain the Hamiltonian from the constrained Lagrangian as follows:  
\bea H_c (Q,S)&=& \left(\langle\dot Q, S\rangle -\La_c (Q,\dot Q)\right)\big|_{\dot Q=
\mathbb{F}\mcal{L}^{-1}(S)} \nn \\ &=&\lin S\Lam^{-1},S\Lam^{-1}\rin- \frac12\lin S
\Lam^{-1},S\Lam^{-1}\rin-\frac12\langle B,QQ\T-I_n\rangle\big|_{\dot Q=S\Lam^{-1}} \nn \\ 
&=&\frac12\lin S\Lam^{-1},S\Lam^{-1}\rin- \frac12\langle B,QQ\T-I_n\rangle\big|_{\dot Q=
S\Lam^{-1}} \nn \\ &=& \frac12\Tr (S\T S\Lam^{-1})- \frac12\Tr(BQQ\T-B). \la{Hamc} \eea 
The extremal (geodesic) flows are then be obtained from this Hamiltonian once the 
flow is restricted to $V(n,N)$, as: 
\bea \dot Q &=& \mbox{Grad}_S H_c(Q,S)= S\Lam^{-1}, \nn \\ 
\dot S &=& -\mbox{Grad}_Q H_c(Q,S)= BQ, \la{conHeq} \eea 
which are identical to equations (\ref{QSeq}) in Proposition \ref{varexs}. Thus, we 
can parametrize $T^\star_Q V(n,N)$ by $(Q,S)$, and equations (\ref{QSeq}) then 
express the extremal flows in terms of this parametrization. \\
 
In the following subsection, we define a momentum quantity that generalizes the 
body momentum of the $N$-dimensional rigid body, and that can be expressed in terms 
of $Q$ and $S$. This expression of the body momentum is then used to solve for 
the unique equivalence class $[U]$, that corresponds to a given extremal solution 
pair $(Q,S)$ of the variational problem. 
%In this subsection, we have expressed the momenta $m$ and $M$ in terms of $Q$ and $S$ 
%and we can also consider these expressions to be the definitions for $m$ and $M$. In a 
%later section, we will give an explicit form for the inverse map, i.e., the $S$ 
%variable in terms of the momenta $m$ and $M$, after we relate the solution to this 
%variational problem in $(Q,S)$ to a solution of the optimal control problem.

\subsection{The Momentum Equation and its Solution}

Using equation (\ref{Qdot}), we can write the Lagrangian in the form
\[ \La(Q,U)= -\frac12 \Tr(\Lam U Q\T QU)= \frac12\langle QU,QU\Lam\rangle. \]
We define the (body) momentum as the $U$-gradient of the Lagrangian with respect to 
the pairing $\langle \cdot,\cdot\rangle$ in (\ref{pair}): 
\be M= 2\mbox{Grad}_U \La(Q,U)= Q\T Q U \Lam+ \Lam U Q\T Q= Q\T S- S\T Q. \la{Mom} \ee 
We call this the body momentum, since in the case $n=N$ ($N$-dimensional rigid 
body), this quantity is the momentum expressed in the body coordinate frame. 
The body and spatial momenta quantities are related by 
\be m= Q M Q\T, \la{Mmrel} \ee
as can be verified using equations (\ref{consm}) and (\ref{Mom}) and the constraint 
$QQ\T=I_n$. One can verify that with the body momentum defined as in (\ref{Mom}), 
equation (\ref{solVP}) is equivalent to 
\be \dot M= [M,U]-A,\;\ \mbox{where } A=Q\T S U+ US\T Q= [Q\T Q, U\Lam U]. 
\la{VPsol} \ee 

We show how these equations generalize the classical rigid body in $N$ dimensions, 
i.e., the case $n=N$. In this case, the Euler-Arnold equations are given by 
\[ \dot Q= QU,\;\ \dot M= [M,U],\ \mbox{ where }\ M=U\Lam+\Lam U, \] 
and $Q\in SO(N)$. From equation (\ref{sVP}), we obtain in the case $n=N$
\[ \dot U\Lam+ \Lam\dot U=  \Lam U^2-U^2\Lam. \]
Taking a time derivative of $M=U\Lam+\Lam U$ and subsituting the above equation, 
we see that $M$ satisfies the Euler-Arnold equation. This is also equivalent to the 
flow of $S$ being represented by (\ref{Sd}). \\

Now we present a solution of the algebraic equation (\ref{Mom}), which we rewrite 
below as
\be M= Q\T S-S\T Q= Q\T QU\Lam+\Lam UQ\T Q\triangleq J_Q(U), \la{JQU} \ee
where $QQ\T= I$, $U=- U\T$, and $\Lam > 0$ (diagonal and positive definite). 
The map $J_Q: \mfrak{so}(N) \rightarrow \mfrak{so}(N)$ is defined by equation 
(\ref{JQU}). We first present a few lemmas, which are necessary to prove 
the main result. 
\begin{lemma}
If $K$ is symmetric and positive definite and $X\in \mfrak{so}(N)$, the map $J: 
\mfrak{so}(n)\rightarrow \mfrak{so}(n)$ given by $J: X \mapsto KX+XK$ has kernel 
zero, and is hence an isomorphism. 
\la{JXlem} 
\end{lemma}
\noindent{\em Proof}: The proof of this statement is similar to the proof of Lemma 
\ref{JBlem}. Since $X$ is skew, there exists a unitary matrix $L$, i.e., 
$L\bar{L}\T=\bar{L}\T L= I_N$, such that $LX\bar{L}\T= \imath \Sigma$, where $\Sigma$ 
is a real diagonal matrix. Thus if $KX+XK=0$, then 
\[ LK\bar{L}\T\Sigma+ \Sigma LK\bar{L}\T=0, \]
where $LK\bar{L}\T$ is a positive definite Hermitian matrix. If $e$ is an eigenvector 
of $\widehat{K}= LK\bar{L}\T$, then 
\[ \widehat{K} e= \lambda e,\quad \lambda>0, \]
and 
\[ \widehat{K}\Sigma+ \Sigma\widehat{K}=0 \Rightarrow \widehat{K}(\Sigma e)+\lambda 
(\Sigma e)= 0. \]
Hence $\Sigma e$ is also an eigenvector of $\widehat{K}$ with eigenvalue $-\lambda 
<0$. But $\widehat{K}$ is positive definite and so all its eigenvalues are strictly 
positive. Thus, we have a contradiction, unless $\Sigma=0$ and hence $X=0$. \qed 

\begin{lemma}
If $J_Q(U)=0$, then $QU=0$.
\la{JQlem}
\end{lemma}
\noindent{\em Proof}: From the proof of Lemma 1, we know that $LU\bar{L}\T=\imath
\Sigma$ where $L$ is unitary. We define the Hermitian matrices  
\[ R=LQ\T Q\bar{L}\T,\quad K=L\Lam\bar{L}\T. \] 
Then the result to be proved becomes equivalent to the following
\be R\Sigma K+ K\Sigma R= L J_Q(U)\bar{L}\T= 0\ \Rightarrow\ R\Sigma=0. \la{lem2} \ee
This is because 
\[ R\Sigma=0\ \Rightarrow\ Q\T QU=0\ \Rightarrow\ QU=0, \]
since $Q\T$ has full (column) rank. Thus we need to show that equation (\ref{lem2}) 
is satisfied, to complete the proof.

Note that $\Sigma$ is real and diagonal, $K$ is positive definite, and $R^2=R$, i.e., 
$R$ is a projection matrix. Since $R$ is a projection onto a subspace of dimension 
$n$, there exists a unitary matrix $T$ ($\bar{T}\T T= T\bar{T}\T= I_N$), such that 
\[ TR\bar{T}\T=  \matl{c|c} \overbrace{I_n}^{n\times n} & 0 \\ \hline 0 & 
\underbrace{0}_{(N-n)\times (N-n)} \matr. \] 
Relative to this decomposition, we write
\beas \lefteqn{T\Sigma\bar{T}\T= \matl{cc} \Sigma_1 & \bar{\Sigma_2}\T \\ \Sigma_2 & 
\Sigma_3 \matr, \;\ \bar{\Sigma_1}\T=\Sigma_1,\;\ \bar{\Sigma_3}\T=\Sigma_3,} \\
\lefteqn{TK\bar{T}\T= \matl{cc} K_1 & \bar{K_2}\T \\ K_2 & K_3 \matr,\;\ \bar{K_1}\T= 
K_1, \;\ \bar{K_3}\T= K_3.} \eeas
We may then express the equation $R\Sigma K+ K\Sigma R= 0$ in the form 
\be \matl{cc} \Sigma_1 & \bar{\Sigma_2}\T \\ 0 & 0 \matr \matl{cc}  K_1 & \bar{K_2}\T 
\\ K_2 & K_3 \matr + \matl{cc} K_1 & \bar{K_2}\T \\ K_2 & K_3 \matr \matl{cc} \Sigma_1 
& 0\\ \Sigma_2 & 0 \matr= 0. \la{RSigS} \ee 
Note that 
\[ R\Sigma=0\ \Leftrightarrow\ \Sigma_1=\Sigma_2=0. \]
Since $K$ is Hermitian and positive definite, 
\[ x\T K_3 x= [0\;\ x\T] \matl{cc} K_1 & \bar{K_2}\T \\ K_2 & K_3 \matr \matl{c} 0 \\ 
x \matr > 0 \]
for all $x\ne 0$. Hence, $K_3$ is Hermitian and positive definite. Also, since $Z=
TK\bar{T}\T$ is positive definite, this implies that $Z^{-1}$ is positive definite, 
and we denote 
\[ Z^{-1}= \matl{cc} L_1 & \bar{L_2}\T \\ L_2 & L_3 \matr, \]
where $L_1$ and $L_3$ are also positive definite (by the above argument). From the 
relations 
\[ K_2 L_1+ K_3 L_2 = 0,\;\ K_1 L_1+\bar{K_2}\T L_2= I_n, \] 
we obtain $L_1$ as 
\[ L_1= (K_1 - \bar{K_2}\T K_3^{-1} K_2)^{-1}. \]
Hence $(K_1-\bar{K_2}\T K_3^{-1}K_2)$ is also positive definite. Now equation 
(\ref{RSigS}) is equivalent to the following two independent equations: 
\beas \Sigma_1 K_1+\bar{\Sigma_2}\T K_2+ K_1\Sigma_1+\bar{K_2}\T\Sigma_2=0, \\
K_2\Sigma_1+ K_3\Sigma_2= 0. \eeas 
We know $K_3$ is invertible so $\Sigma_2= -K_3^{-1} K_2\Sigma_1$ and 
\beas \Sigma_1 K_1+ (-\overline{K_3^{-1} K_2\Sigma_1})\T K_2+K_1\Sigma_1+ \bar{K_2}\T 
(-K_3^{-1}K_2\Sigma_1)&=&0 \\ 
\Rightarrow \Sigma_1 (K_1- \bar{K_2}\T K_3^{-1}K_2)+ (K_1- \bar{K_2}\T K_3^{-1}K_2)
\Sigma_1 &=& 0. \la{Lyte} \eeas 
Here we note that $K_1- \bar{K_2}\T K_3^{-1}K_2$ is a Hermitian positive definite 
matrix. Using the result of Lemma \ref{JXlem}, we see that equation (\ref{Lyte}) implies 
that $\Sigma_1=0$, and hence $\Sigma_2= -K_3^{-1} K_2\Sigma_1= 0$. Thus $R\Sigma=0$, 
which implies as we have shown that $QU=0$. \qed 

Now we present the solution to the algebraic equation (\ref{Mom}) or (\ref{JQU}). 
\begin{theorem}
All solutions of the equation (\ref{JQU}) have the form 
\be U=U_1+U_2+ V, \la{Usol} \ee
where 
\bea U_1&=& \Lam^{-1}\big(Q\T(Q\Lam^{-1}Q\T)^{-1}S -S\T (Q\Lam^{-1}Q\T)^{-1}Q\big)
\Lam^{-1}, \nn \\ 
U_2 &=& \Lam^{-1}(Q\T \wh{J_Q}^{-1}(R) Q)\Lam^{-1}, \la{U1U2} \eea
and 
\bea R&=& (Q\Lam^{-1}S\T)(Q\Lam^{-1}Q\T)^{-1}-(Q\Lam^{-1}Q\T)^{-1}(S\Lam^{-1}Q\T), \nn 
\\ \wh{J_Q}(X)&=& (Q\Lam^{-1}Q\T) X+ X(Q\Lam^{-1} Q\T);\;\ \wh{J_Q}: \mfrak{so}(n)
\rightarrow \mfrak{so}(n) \la{Req} \eea
and $QV=0$. 
\la{algeq}
\end{theorem} 
\noindent{\em Proof}: We evaluate $J_Q(U_1)$ as follows 
\beas \fl J_Q(U_1)&= Q\T Q\Lam^{-1}\big( Q\T(Q\Lam^{-1}Q\T)^{-1}S -S\T (Q\Lam^{-1}
Q\T)^{-1}Q\big)+ \big(Q\T(Q\Lam^{-1}Q\T)^{-1}S \\  & -S\T (Q\Lam^{-1}Q\T)^{-1}Q
\big) \\ &= Q\T S- Q\T(Q\Lam^{-1}S\T)(Q\Lam^{-1}Q\T)^{-1}Q+ Q\T(Q\Lam^{-1}Q\T)^{-1}
(S\Lam^{-1}Q\T)Q-S\T Q \\ &= M-Q\T RQ. \eeas
Let $U_2=\Lam^{-1}Q\T XQ\Lam^{-1}$; then  
\[ J_Q(U_2)= Q\T(Q\Lam^{-1}Q\T)XQ+ Q\T X(Q\Lam^{-1}Q\T)Q= Q\T\wh{J_Q}(X)Q. \]
Thus if $\wh{J_Q}(X)=R$, then $J_Q(U_2)= Q\T RQ$. From Lemma \ref{JXlem}, we know 
that $\wh{J_Q}: \mfrak{so}(n)\rightarrow \mfrak{so}(n)$ is an isomorphism. It follows 
that $\wh{J_Q}$ is invertible, $\wh{J_Q}^{-1}(R)=X$ is unique, and $U_2=\Lam^{-1}Q\T 
XQ\Lam^{-1}$ satisfies $J_Q(U_2)= Q\T RQ$. Thus, 
\[ J_Q(U_1+U_2)=J_Q(U_1)+J_Q(U_2)= M-Q\T RQ+ Q\T RQ= M, \]
and hence all solutions of (\ref{JQU}) have the form
\[ U= U_1+U_2+V,\;\ J_Q(V)= 0. \]
From Lemma \ref{JQlem}, we know that $J_Q(V)=0 \Rightarrow QV=0$. Thus, we have 
proved this theorem. \qed \\
Note that this solution to the algebraic equation (\ref{Mom}) or (\ref{JQU}) uses 
the decomposition of $M$ in the $(Q,S)$ variables. \\

%\subsection{The Extremal Flows in terms of $Q$ and $S$} 
%Note that although the solution of the algebraic equation (\ref{Mom}) or (\ref{JQU}) 
%is unique only up to the equivalence class $[U]\subset \mfrak{so}(N)$, this determines 
%$S=QU\Lam$ uniquely, since $QU$ is uniquely determined from this solution. As we will 
%see in the next section, the body momentum $M$ can be decomposed in terms of the 
%extremal solutions to the optimal control problem in exactly the same manner as 
%equation (\ref{Mom}), with $S$ replaced by the extremal costate $P$. Thus, the 
%solution $U$ to this momentum equation can also be expressed in terms of the extremal 
%solution to the optimal control problem. Now we express the flow of the variable $S$ 
%along the extremal solutions to the variational problem, given the flow in terms of 
%$Q$ and $U$. \\
Let $X=\wh{J_Q}^{-1}(R)$, where $R$ and $\wh{J_Q}$ are as defined in equation 
(\ref{Req}). The following lemma gives an expression for $X$ in terms of $Q$ and $S$. 
\begin{lemma}
If $\wh{J_Q}$ and $R$ are as given by equation (\ref{Req}), then we can express 
\be X=\wh{J_Q}^{-1}(R)=(Q\Lam^{-1}Q\T)^{-1} Q\Lam^{-1}S\T (Q\Lam^{-1}Q\T)^{-1}. 
\la{Xeq} \ee
\la{Xlem} 
\end{lemma} 
\vspace*{-2mm}
\noindent {\em Proof}: If $U$ is a solution of the algebraic equation (\ref{JQU}), 
then the quantity $S$ can be expressed in terms of $U$ of as 
\beas \fl S= QU\Lam = Q\Lam^{-1}\{Q\T (Q\Lam^{-1}Q\T)^{-1}S-S\T
(Q\Lam^{-1}Q\T)^{-1}Q\}+ Q\Lam^{-1} Q\T XQ \\
\Rightarrow S= S-Q\Lam^{-1}S\T (Q\Lam^{-1}Q\T)^{-1}Q+ Q\Lam^{-1}Q\T XQ \\ 
\Rightarrow Q\Lam^{-1}Q\T XQ = Q\Lam^{-1}S\T (Q\Lam^{-1}Q\T)^{-1}Q \\ 
\Rightarrow XQ= (Q\Lam^{-1}Q\T)^{-1} Q\Lam^{-1}S\T (Q\Lam^{-1}Q\T)^{-1}Q. \eeas 
Post-multiplying both sides of the above expression with $Q\T$, we get the 
expression in (\ref{Xeq}) for $X$. \qed \\
One can verify that $\wh{J_Q}(X)=R$, where $X$  is as given by (\ref{Xeq}), 
as follows
\beas \wh{J_Q}(X)&=&(Q\Lam^{-1}Q\T)X+X(Q\Lam^{-1}Q\T) \\
&=& (Q\Lam^{-1}S\T)(Q\Lam^{-1}Q\T)^{-1}+(Q\Lam^{-1}Q\T)^{-1}(Q\Lam^{-1}S\T) \\
&=& (Q\Lam^{-1}S\T)(Q\Lam^{-1}Q\T)^{-1}-(Q\Lam^{-1}Q\T)^{-1}(S\Lam^{-1}Q\T) 
=R, \eeas
where we used $Q\Lam^{-1}S\T= -QUQ\T$ is skew-symmetric in the last step above. \\

%Note that the relation $\wt{J_Q}(B)= -2\dot Q\dot Q\T$ is like a Lyapunov equation, 
%where the right-hand side is a symmetric, negative definite matrix. However, in this 
%case, the solution we desire is the matrix $B$, which is also symmetric. 
The following statement is a corollary of Theorem \ref{algeq} and Lemma \ref{Xlem}.
\begin{corollary}
For a given $Q$, the map $\mcal{Z}: T^\star_Q V(n,N)\rightarrow T_Q V(n,N)$ where 
$\mcal{Z}: S\mapsto [U]$ is given by Theorem \ref{algeq}, is an isomorphism. 
\la{S2Uc}
\end{corollary}
\noindent {\em Proof}: Clearly, the map $\mcal{Z}^{-1}: [U]\mapsto S=QU\Lam$ is a 
linear isomorphism. Observe from Theorem \ref{algeq} and Lemma \ref{Xlem}, that we 
get the following expressions for $U_1$ and $U_2$:
\beas U_1 &=& \Lam^{-1}\left[ Q\T(Q\Lam^{-1}Q\T)^{-1}S- S\T(Q\Lam^{-1}Q\T)^{-1}
Q \right]\Lam^{-1}, \\
U_2&=& \Lam^{-1}\left[Q\T(Q\Lam^{-1}Q\T)^{-1}Q\Lam^{-1}S\T(Q\Lam^{-1}Q\T)^{-1}Q
\right]\Lam^{-1}, \eeas
which depend linearly on $S$. 
%\beas U_1+U_2 &=& \Lam^{-1}\left[ Q\T(Q\Lam^{-1}Q\T)^{-1}S- S\T(Q\Lam^{-1}Q\T)^{-1}
%Q \right. \\ & &\left. + Q\T(Q\Lam^{-1}Q\T)^{-1}Q\Lam^{-1}S\T(Q\Lam^{-1}Q\T)^{-1}Q
%\right]\Lam^{-1}, \eeas
It is now easy to verify that 
\beas QU_1 &=& [S- Q\Lam^{-1}S\T(Q\Lam^{-1}Q\T)^{-1}Q]\Lam^{-1}, \\
QU_2 &=& [Q\Lam^{-1}S\T(Q\Lam^{-1}Q\T)^{-1}Q]\Lam^{-1}, \eeas 
which gives us $Q(U_1+U_2)=S\Lam^{-1}$, validating the relation $S=QU\Lam$. 
Hence, the map $\mcal{Z}: S\mapsto [U]$ is an isomorphism. \qed \\
Thus, Theorem \ref{algeq} and Lemma \ref{Xlem} describe the exact relationship 
between the extremal solutions expressed in terms of $(Q,[U])$ in equations 
(\ref{Qdot})-(\ref{solVP}), and those expressed in terms of $(Q,S)$ in equations 
(\ref{QSeq})-(\ref{JQB}). 

\subsection{The Discrete Variational Problem}

The discrete counterpart of the variational problem (\ref{VPnN}) and the discrete 
extremal trajectories obtained thereof were given by Moser and Veselov (1991). 
The discrete variational problem is given by 
\be \min_{Q_k}\sum_k\frac12\langle Q_{k+1}\Lam,Q_k \rangle, 
\label{dVPnN} \ee
subject to $Q_k Q_k\T= I_n$. The extremal trajectories to this discrete variational 
problem are given by Moser and Veselov (1991)
\be Q_{k+1}\Lam+ Q_{k-1}\Lam= B_k Q_k,\;\ k\in\mathbb{Z} \la{disvar} \ee
where $B_k= B_k\T \in\bbr^{nn}$ is a (symmetric) Lagrange multiplier matrix for the 
symmetric constraint $Q_k Q_k\T=I_n$. The above equation is the discrete counterpart 
of equation (\ref{QSeq}). The discrete body momentum is defined as 
\be M_k= Q_{k-1}\T Q_k\Lam-\Lam Q_k\T Q_{k-1}. \la{dbMom} \ee 
Since $B_k$ is symmetric, equation (\ref{disvar}) is equivalent to the conservation 
of the discrete spatial momentum 
\be m_{k+1}\triangleq Q_{k+1}\Lam Q_k\T- Q_k\Lam Q_{k+1}\T= Q_k\Lam Q_{k-1}\T-
Q_{k-1}\Lam Q_k\T \triangleq m_k. \la{dsmom} \ee
Thus, $m_k= Q_{k-1}M_k Q_{k-1}\T$ is conserved along the discrete extremal 
trajectories. Let us define 
\[ U_k\triangleq Q_{k-1}\T Q_k, \] 
which implies that 
\[ Q_k= Q_{k-1}U_k,\;\ M_k= U_k\Lam- \Lam U_k\T. \]
The following proposition gives the discrete extremal trajectories in terms of $U_k$ 
and the discrete body momentum $M_k$. 
\begin{proposition}
The extremal trajectories of the discrete variational problem (\ref{dVPnN}) on the 
Stiefel manifold $V(n,N)$ in terms of $(M_k, U_k)$ are given by: 
\be \lefteqn{M_{k+1}=U_k\T M_k U_k- A_k,} \la{dsolVP} \ee 
where 
\be A_k= U_k\T \Lam\big(I_N-U_k\T U_k\big) -\big(I_N-U_k\T U_k\big)\Lam U_k. 
\la{Akeq} \ee
\la{discVP}
\end{proposition}
\noindent{\em Proof:} One can obtain this result from the second order difference 
equation (\ref{disvar}) which gives the extremal trajectories for (\ref{dVPnN}). 
Using this equation, we can represent the body momentum at the $(k+1)$-th step as 
\bea M_{k+1}&=& Q_k\T (B_k Q_k-Q_{k-1}\Lam)-(B_k Q_k-Q_{k-1}\Lam)\T Q_k \nn \\
&=& \Lam U_k- U_k\T \Lam. \la{dbmkp1} \eea
From equation (\ref{dsmom}) expressing conservation of the spatial momentum, we get 
\beas Q_k M_{k+1}Q_k\T &=& Q_{k-1}M_k Q_{k-1}\T= Q_k U_k\T M_k U_k Q_k\T \\ 
\Leftrightarrow M_{k+1}&=& U_k\T M_k U_k- A_k\;\ 
\mbox{where }\ Q_k A_k Q_k\T= 0. \eeas 
From the above expression, we get 
\[ A_k= M_{k+1}-U_k\T M_kU_k, \] 
and now using equation (\ref{dbmkp1}), we obtain $A_k$ as given in equation 
(\ref{Akeq}). This proves the given result. \qed \\
Note that equation (\ref{dbmkp1}) is equivalent to equations (\ref{dsolVP})-(\ref{Akeq}) 
for the discrete body momentum. These equations are therefore the discrete counterpart 
of equation (\ref{VPsol}) for the continuous case. \\ %Thus, Proposition \ref{discVP} 
%expresses the discrete extremal flow in $T^\star V(n,N)$ parametrized locally by 
%$(Q_k,Q_k M_k)$. \\

Proposition \ref{discVP} can also be used to prove Theorem 4 of Moser and Veselov 
(1991), which gives the following set of isospectral deformations for the discrete 
extremal flows 
\be L_{k+1}(\lambda)= C_k(\lambda)L_k(\lambda)C_k^{-1}(\lambda), \la{isosp} \ee
where $L_k(\lambda)=\Lam^2+\lambda M_k-\lambda^2 Q_{k-1}\T Q_{k-1}$, and $C_k(\lambda)
=\Lam-\lambda Q_k\T Q_{k-1}$. Note that $L_k(\lambda)$ can be factored as 
\[ L_k(\lambda)= (\Lam+\lambda Q_{k-1}\T Q_k)(\Lam-\lambda Q_k\T Q_{k-1})= 
C_k\T (-\lambda)C_k(\lambda). \] 
Hence, the determinant of $L_k$ is an even polynomial in $\lambda$ of degree $2n$.  
Using equations (\ref{dsolVP})-(\ref{Akeq}), we can express the left hand side of 
equation (\ref{isosp}) as 
\[ \fl \Lam^2+\lambda (Q_k\T Q_{k-1}M_k Q_{k-1}\T Q_k- A_k)-\lambda^2 Q_k\T Q_k= 
\Lam^2+\lambda (\Lam Q_{k-1}\T Q_k- Q_k\T Q_{k-1}\Lam)-\lambda^2 Q_k\T Q_k. \] 
The right hand side of equation (\ref{isosp}) is obtained from the factorization 
of $L_k(\lambda)$ as follows
\beas C_k(\lambda)\big(C_k\T (-\lambda)C_k (\lambda)\big)C_k^{-1}(\lambda)&=& 
C_k(\lambda) C_k\T (-\lambda) \\
&=& (\Lam-\lambda Q_k\T Q_{k-1})(\Lam+\lambda Q_{k-1}\T Q_k) \\
&=& \Lam^2+ \lambda (\Lam Q_{k-1}\T Q_k- Q_k\T Q_{k-1}\Lam)-\lambda^2 Q_k\T Q_k. \eeas 
Using this set of isospectral deformations, Moser and Veslov (1991) 
give a method to reconstruct the discrete flow under further conditions, and 
the discrete flow is integrable in this sense. 

\section{The Optimal Control Problem}\la{sec:ocpro}

We now study the Hamiltonian approach to the variational problem (\ref{OCPnN}). 
The Hamiltonian for the optimal control problem (\ref{OCPnN}) is given by 
\bea H(P,Q,U)&=& \langle P,QU \rangle- \La(Q,U) \nn \\
&=& \Tr (P\T QU)+\frac12 \Tr (\Lam U Q\T QU), \la{Ham} \eea 
where $P\in \bbr^{nN}$ denotes the costates (Lagrange multipliers). This 
optimal control problem is nominally posed on $W_{n,N}=\bbr^{nN}\times \bbr^{nN}$, 
on which the symplectic structure is given by the symplectic form
\be \omega((A_1,A_2),(B_1,B_2))=\langle A_1,B_2\rangle-\langle A_2,B_1\rangle. 
\la{omWnN} \ee 
We restrict the solutions of this optimal control problem to those which are 
governed by extremals that leave the submanifolds $W^k_{n,N}$ invariant, where 
the $W^k_{n,N}$ are level sets of $W_{n,N}$ specified by 
\be W^k_{n,N}= \{(Q,P)\in W_{n,N};\;\ QQ\T= I_n;\;\ QP\T+PQ\T= k\}, \la{WknN} \ee
where $k$ is some constant symmetric $n\times n$ matrix. %When we choose a particular 
%solution pair $(Q,P)$, we are also implicitly choosing the solution submanifold 
%$W^k_{n,N}$. 
Note that the dimension of $W^k_{n,N}$ is given by 
\be \mbox{Dim }W^k_{n,N}= \mbox{Dim }V(n,N)+nN-\frac{n(n+1)}{2}= 2\mbox{Dim }V(n,N), 
\la{dimWk} \ee
which is equal to the dimension of the (co)tangent bundle of $V(n,N)$. 

\subsection{Space of extremal solutions to the optimal control problem}

%The space of solutions of the optimal control problem is given by $W^k_{n,N}$, 
%as defined by (\ref{WknN}), where $k$ is a fixed symmetric $n\times n$ matrix. 
%We showed in Section 2.2, equation (\ref{dimWk}), that the dimension of this 
%space is twice that of the Stiefel manifold. 
Consider the vector space $L$ of vector fields on $W_{n,N}$ characterized by the 
differential equations 
\bea \dot Q&=& QU, \;\ U\in\mfrak{so}(N), \nn \\
\dot P&=& PU+QV, \;\ V\in\mfrak{so}(N). \la{vecf} \eea  
Let $X_1=(QU_1,PU_1+QV_1)$ and $X_2=(QU_2,PU_2+QV_2)$ be vectors in $L$. 
This vector space is seen to be a Lie algebra, since:
\beas [X_1,X_2] &=& [(QU_1, PU_1+QV_1), (QU_2,PU_2+QV_2)]\\ %=\frac{d}{ds_1}
%\frac{d}{ds_2} g(s_1)h(s_2)g(s_1)^{-1} \big|_{s_1=0,s_2=0}  
&=& (Q[U_1,U_2],P[U_1,U_2]+ Q[V_1,U_2]+Q[U_1,V_2]), \eeas 
also belongs to $L$. We now show that the submanifolds $W^k_{n,N}$ are integral 
manifolds to the involutive distribution of vector fields in the Lie algebra $L$.
\begin{lemma}
The submanifolds $W^k_{n,N}\subset W_{n,N}$ are integral submanifolds to the 
involutive distribution on $W_{n,N}$ defined by the Lie algebra of vector 
fields, $L$. 
\la{LiWknN}
\end{lemma} 
\noindent {\em Proof:} Differentiating the constraints $QQ\T=I_n$ and $PQ\T+QP\T=k$, 
along trajectories of vector fields in $L$ defined by the system (\ref{vecf}), 
we get: 
\beas QUQ\T- QUQ\T=0, \\
(PU+QV)Q\T- PUQ\T+ QUP\T+ Q(-UP\T-VQ\T)=0. \eeas 
This shows that all vector fields in $L$ are tangent to each of the submanifolds 
$W^k_{n,N}$ of $W_{n,N}=\bbr^{nN}\times\bbr^{nN}$. We find the dimension of 
the subspace $L(Q,P)$ of $T_{(Q,P)} W^k_{n,N}$ spanned by the vector fields in 
$L$. Since $U$ and $V$ are independent 
\beas & &\mdim \{(QU,PU+QV),\ (Q,P)\in W^k_{n,N},\ U,V\in\mfrak{so}(N) \} \\
&=& \mdim \{(QU,PU),\ U\in\mfrak{so}(N)\}+ \mdim\{ (0,QV),\ V\in\mfrak{so}(N)\} \\
&=& \mdim \{(QU,0),\ U\in\mfrak{so}(N)\}+ \mdim\{ (0,QV),\ V\in\mfrak{so}(N)\} \\
&=& 2\mdim \{QU,\ U\in\mfrak{so}(N),\ QQ\T=I_n\} = 2\mdim V(n,N), \eeas
since the dimension of the tangent space to $V(n,N)$ at $Q$ is the same as the 
dimension of $V(n,N)$ itself. %, we get 
%\[ \mdim L(Q,P)= 2\mdim V(n,N). \]
From (\ref{dimWk}), we see that this is also the dimension of $W^k_{n,N}$, and 
hence that of $T_{(Q,P)} W^k_{n,N}$. Thus, we conclude that at every $(Q,P)\in 
W^k_{n,N}$, 
\[ L(Q,P)= T_{(Q,P)} W^k_{n,N}. \]
This proves the given result. \qed \\

Now that we have shown that the $W^k_{n,N}$ are integral submanifolds of the Lie 
algebra $L$, we next state and prove the following important result.
\begin{theorem}
The space $W^k_{n,N}$ defined by (\ref{WknN}) is a symplectic submanifold 
of $W_{n,N}=(\bbr^{nN}\times \bbr^{nN},\omega)$.
\la{sympman}
\end{theorem}
\noindent {\em Proof:} To show this result, we need to show that the symplectic form 
$\omega$ on $W_{n,N}$ given by (\ref{omWnN}), is non-degenerate on $W^k_{n,N}$ 
when restricted to $W^k_{n,N}$. Let us denote the restriction of the symplectic 
form on $W_{n,N}$ to $W^k_{n,N}$ by 
\[ \Omega= \omega |_{W^k_{n,N}}. \]
Tangent vectors to $W^k_{n,N}$ are given by (\ref{vecf}). Thus we must show that 
\[ \Omega ((QU_1,PU_1+QV_1),(QU_2,PU_2+QV_2))=0 \]
for all $(U_2,V_2)\in \mfrak{so}(N)\times\mfrak{so}(N)$ implies that $(QU_1,
PU_1+QV_1)=(0,0)$. By definition of $\Omega$, we have 
\be \langle QU_1,PU_2+QV_2\rangle -\langle QU_2,PU_1+QV_1\rangle= 0,\ \forall\ 
(U_2,V_2)\in \mfrak{so}(N)\times\mfrak{so}(N). \la{degsf} \ee 
Setting $U_2=0$, we get 
\[ \Tr (U_1 Q\T QV_2)=0,\ \forall\ V_2. \] 
Since $V_2$ is skew-symmetric, this in turn implies that $U_1Q\T Q$ is symmetric, 
or \[ QU_1 Q\T Q= -QU_1,\ \mbox{ since } QQ\T=I_n. \]
Thus, each row of $QU_1$ is an eigenvector of $Q\T Q$ with an eigenvalue of $-1$. 
Since $Q\T Q$ is positive semi-definite, this is a contradiction unless $QU_1=0$. 
Now setting $QU_1=0$ we have from (\ref{degsf}) that 
\[ \langle QU_2,PU_1+QV_1\rangle=\Tr (U_2Q\T (PU_1+QV_1))=0,\ \forall\ 
U_2. \]
Since $U_2$ is skew-symmetric, this implies that 
\beas Q\T (PU_1+QV_1)&=& -(U_1 P\T+ V_1 Q\T)Q \\
\Rightarrow PU_1 Q\T+ QV_1 Q\T&=& -QU_1 P\T- QV_1 Q\T \\
\Rightarrow QV_1 Q\T &=& -QV_1 Q\T\ \mbox{ since } QU_1=0 \\
\Rightarrow QV_1 Q\T &=& 0 \\
\Rightarrow Q\T (PU_1+QV_1)Q\T &=& 0= -(U_1P\T+V_1Q\T)QQ\T \\ 
\Rightarrow PU_1+ QV_1 &=& 0. \eeas 
Thus, $(\Omega,W^k_{n,N})$ for any value of $k$ is a symplectic submanifold 
of $(\omega, W_{n,N})$. \qed \\
We restrict the extremal flows of the optimal control problem to the symplectic 
manifolds $W^k_{n,N}$, for a value of $k$ given by the initial conditions.
%In the following subsection, we give the extremal solutions to the optimal control 
%problem in $W^k_{n,N}$ and, later on, show that $W^0_{n,N}$ (where $k=0$) is 
%symplectomorphic to the cotangent bundle, $T^\star V_{n,N}$. 

\subsection{Solution to the Optimal Control Problem}

The extremal solutions on $W^k_{n,N}$ to the optimal control problem are 
characterized by the following result. 
\begin{theorem}
The extremal trajectories of the optimal control problem (\ref{OCPnN}), $(Q,P)\in 
W^k_{n,N}$, are given by 
\bea \lefteqn{\dot Q= QU, \;\ U\in \mfrak{so}(N),} \la{sOCP} \\ 
\lefteqn{\dot P= PU-QA,} \la{solOCP} \eea %= PV,} 
where $A= Q\T Q U\Lam U- U\Lam U Q\T Q$. 
\la{exsolW}
\end{theorem}
\noindent {\em Proof}: The Hamiltonian (\ref{Ham}) can be written in the alternate 
form
\be H(P,Q,U)= \langle P,QU\rangle+ \frac12\langle Q\T Q, U\Lam U\rangle %\nn \\
%&=&\frac12 \langle Q\T P-P\T Q, U\rangle+ \frac12 \langle Q\T Q, U\Lam U\rangle. 
\la{Hama} \ee
This optimal control problem may be restricted (with possible loss of generality) 
so that extremal trajectories lie on the symplectic manifold $(W^k_{n,N},\Omega)$, 
where $\Omega= \omega |_{W^k_{n,N}}$ is the symplectic two-form given by
\be \Omega((\dot Q_1,\dot P_1),(\dot Q_2,\dot P_2))=\langle \dot Q_1,\dot P_2\rangle
-\langle \dot Q_2,\dot P_1\rangle. \la{2fWnN} \ee
Hence, we have
\be dH(X_2)= \Omega(X_1,X_2), \la{dHom} \ee
where $X_1=(\dot Q_1,\dot P_1)$ and $X_2=(\dot Q_2,\dot P_2)$ are vector fields  
such that $\dot Q_i= QU_i$, $\dot P_i= PU_i-QV_i$ $i=1,2$, $U_i,V_i\in so(N)$, and 
$X_1=X_H$. Thus, from (\ref{dHom}), we get 
\beas & & \langle PU_2-QV_2, QU_1\rangle -\langle PU_1-QV_1, QU_2\rangle \\
&=& \langle PU_2-QV_2, QU\rangle +\langle P,QU_2 U\rangle+ \frac12\langle 
[Q\T Q,U_2],U\Lambda U\rangle. \eeas
Setting $U_2=0$ in the above expression, we get $QU_1=QU$. Thus, from above, we get 
\beas & & \langle PU_2,QU\rangle-\langle PU-QV_1,QU_2\rangle \\
&=& \langle PU_2,QU\rangle+\langle P, QU_2 U\rangle+\frac12\langle [Q\T Q,U_2],
U\Lam U\rangle. \eeas
So 
\beas \langle QV_1,QU_2\rangle&=& \frac12\langle [Q\T Q,U_2],U\Lam U\rangle \\
&=& \frac12\langle Q\T QU_2-U_2 Q\T Q,U\Lam U\rangle \\
&=& \frac12\langle QU_2,QU\Lam U\rangle +\frac12\langle QU_2,QU\Lam U\rangle \\
&=& \langle QU_2,QU\Lam U\rangle. \eeas
But 
\[ \langle QU_2, QU\Lam U\rangle= \langle QU_2, Q(Q\T QU\Lam U-U\Lam UQ\T Q)\rangle, \]
since $QQ\T=I_n$ and 
\[ \langle QU_2, QU\Lam UQ\T Q\rangle= \langle QU_2 Q\T, QU\Lam UQ\T\rangle=0, \]
as $U_2$ is skew-symmetric and $QU\Lam UQ\T$ is symmetric. Hence, we get 
$V_1=[Q\T Q,U\Lam U]$ up to equivalence, and $X_H$ is the vector field given by 
\[ \dot Q= QU,\;\ \dot P=PU-QA,\ \mbox{ where }\ A=[Q\T Q,U\Lam U]. \]
Thus, the Hamiltonian vector field $X_H$ prescribes the flow given by equations 
(\ref{sOCP})-(\ref{solOCP}).  \qed \\

We determine the optimal control $U$ applying Pontryagin's maximum principle 
(see Bloch \etal (2003), Gelfand and Fomin (2000) and Kirk (2004)). The Hamiltonian 
in (\ref{Ham}) can also be expressed as 
\be H(Q,P,U)= \Tr \frac12(P\T Q-Q\T P)U+\frac12\Tr Q\T QU\Lam U. 
\la{Hamm} \ee
Then $\mbox{Grad}_UH(Q,P,U^*)=0$ with respect to the pairing in (\ref{pair}) is 
equivalent to
\be Q\T P-P\T Q=Q\T QU^*\Lam+\Lam U^*Q\T Q=M, \la{maxprin} \ee 
where $U^*$ is the optimal control. This equation gives the momentum $M$ in terms of 
the states $Q$ and costates $P$ and also in terms of $Q$ and $U^*$ (or $Q$ and $S$), 
as given by the solution to the variational problem in equation (\ref{Mom}). We 
now appeal to Theorem \ref{algeq} to give an explicit representation for $U^*$ in 
terms of $(Q,P)$, by replacing $S$ with $P$. \\

Taking the time derivative of the equation $M= Q\T P-P\T Q$ along the vector field 
given by equations (\ref{sOCP})-(\ref{solOCP}), we get 
\be \dot M= [M, U]-A, \la{Madot} \ee
%where $V\in [U]$, 
which is identical to equation (\ref{VPsol}) obtained from the solution to the 
variational problem. \\ %since $Q\T QA+AQ\T Q= A$. 

Since $QAQ\T=0$, one can right multiply (\ref{solOCP}) with $Q\T$ to get 
\[ \dot P Q\T= PUQ\T, \] %= PVQ\T, \]
which implies (from \ref{sOCP}) that $\lambda=PQ\T$ is conserved.  
Hence the symmetric quantity
\be k= PQ\T+ QP\T= \lambda+\lambda\T, \la{keq} \ee
is conserved along the trajectories of (\ref{sOCP})-(\ref{solOCP}).
Since the spatial momentum is
\be m=QMQ\T= PQ\T- QP\T= \lambda-\lambda\T, \la{momm} \ee
as originally defined in (\ref{consm}), this quantity is also conserved along the 
flow of the extremal solution (\ref{sOCP})-(\ref{solOCP}) to the optimal control 
problem. Hence, for each initial condition set, the solution trajectory is 
confined to the level set 
\be \fl W^k_m= \{(Q,P)\in W_{n,N};\;\ QQ\T= I_n;\;\ QP\T+PQ\T= k,\;\ PQ\T-QP\T=m\}, 
\la{Wkm} \ee
which is a submanifold of $W^k_{n,N}$. The manifold $W^k_{n,N}$ is formed by an 
union of the $W^k_m$ over all values of $m$. Clearly, the Hamiltonian vector fields 
given by (\ref{sOCP})-(\ref{solOCP}) are tangent to the $W^k_m$ which are level 
momentum sets of this Hamiltonian flow. 
Note that $PQ\T=\lambda= \frac12 (k+m)$ is constant for $(P,Q)\in W^k_m\subset 
W^k_{n,N}$. \\  

The quantity that corresponds to the body momentum in the full-ranked case, is 
according to (\ref{maxprin})
\be M= Q\T P- P\T Q, \la{bMom} \ee
which has been previously defined in equation (\ref{Mom}) in terms of the solution to 
the variational problem. In the case $n=N$ we have obtained Bloc \etal (2002) the 
following relation when $\|m\| < 2$
\[ P= Q e^{\sinh^{-1}\frac{M}{2}}= e^{\sinh^{-1}\frac{m}{2}} Q. \] 
For the general case, $1\le n< N$, when $PQ\T=\lambda$ is orthogonal and $\|m\|<2$, 
the extremal solutions (\ref{sOCP})-(\ref{solOCP}) in $W^k_m\subset W^k_{n,N}$ can 
be expressed as  
\[ Q\in V(n,N),\;\ P= e^{\sinh^{-1}\frac{m}{2}} Q+ QR, \]
where $k=c\T+c$, $c=e^{-\sinh^{-1}\frac{m}{2}}$, and
\[ R\in \mfrak{so}(N),\;\ QRQ\T=0,\;\ \dot R= [R,U]-A. \]
This can be easily verified by direct substitution into equations 
(\ref{sOCP})-(\ref{solOCP}). \\  

Note that the quantity $k$ specifies the symplectic submanifold on which the extremal 
solution lies, while the spatial momentum $m$ specifies the momentum level set. We 
want to express the costates $P$ of this optimal control problem in terms of $Q$ and 
the momentum quantities $m$ and $M$, given by (\ref{momm}) and (\ref{bMom}) 
respectively. %This choice of the optimal costate $P$ is given by the 
%result in the following discussion. 
Note also that if $(Q,P)$ is an extremal trajectory to the optimal control 
problem, i.e., 
\be \dot Q= QU,\;\ \dot P= PU-QA, \la{QPdot} \ee 
then $(Q,P+hQ)$, where $h$ is a constant $n\times n$ matrix, is also an 
extremal trajectory satisfying (\ref{QPdot}). If $(Q,P)$ lies in the solution 
submanifold $W^k_{n,N}$ given by (\ref{WknN}), then $(Q,P+hQ)\in W^{k+b}_{n,N}$ 
where $b=h+h\T$. Further, if $h$ is symmetric and $(Q,P)$ lies in the momentum 
level set $W^k_m\subset W^k_{n,N}$, then $(Q,P+hQ)$ lies in the momentum level 
set $W^{k+2h}_m\subset W^{k+2h}_{n,N}$; i.e., the spatial momentum value remains 
unchanged. One can verify that in this case, the body momentum value $M$ also 
remains unchanged. %Hence, we have $(Q,P-\frac12 kQ)\in W^0_{n,N}$. 
This brings us to the following remarkable result. 
\begin{proposition}
We have a map $\Xi: W^k_{n,N}\rightarrow W^0_{n,N}$, defined by $\Xi: (Q,P)
\mapsto (Q,P_0)\in W^0_{n,N}$, where 
\be P_0= P-\frac12 kQ, \;\ (Q,P)\in W^k_{n,N}. \la{P0eq} \ee 
Further, if $(Q,P)\in W^k_m\subset W^k_{n,N}$, then $(Q,P_0)\in W^0_m\subset 
W^0_{n,N}$, i.e., it satisfies equations (\ref{QPdot}) and 
\beas P_0 Q\T- QP_0\T = & m &= PQ\T-QP\T,\\ Q\T P_0-P_0\T Q =&M &= 
Q\T P-P\T Q. \eeas
Hence the map $\Xi$ leaves the spatial and body momenta unchanged. 
\la{Wkto0}
\end{proposition} 
\noindent {\em Proof}: Let $(Q,P)\in W^k_{n,N}$. Then we have 
\beas P_0Q\T+ QP_0\T &=& (P-\frac12 kQ)Q\T +Q(P-\frac12 kQ)\T \\
&=& PQ\T -\frac12 k+ QP\T-\frac12 k \\
&=& k-k=0. \eeas
Thus $(Q,P)$ is mapped to $(Q,P_0)\in W^0_{n,N}$ by $\Xi$. If $(Q,P)\in W^k_m$, then 
$PQ\T-QP\T=m$ in addition to $PQ\T+QP\T=k$, and we have 
\beas P_0Q\T-QP_0\T &=& (P-\frac12 kQ)Q\T -Q(P-\frac12 kQ)\T \\
&=& PQ\T -\frac12 k- QP\T+\frac12 k \\ 
&=& PQ\T- QP\T=m. \eeas
This proves the second part of the statement. \qed \\

The following result is a corollary of Proposition \ref{Wkto0}, and gives the 
costate $P=P_0$ as a function of $Q$, $m$ and $M$, such that $(Q,P_0)\in W^0_m$. 
As we show later, this costate is a natural choice since a direct relation exists 
between the symplectic forms on $T^\star V(n,N)$ and $W^0_{n,N}$. 
\begin{corollary}
The pair $(Q,P_0)\in W^0_m$ satisfying equations (\ref{QPdot}), may be expressed as 
\be P_0= -\frac1{2} mQ+ QM. \la{Peqn} \ee
%Then
%\[ m=P_0 Q\T- QP_0\T,\;\ M=Q\T P_0- P_0\T Q. \] 
For this solution, $P_0=Q\bar{M}$ where 
\be \bar{M}= -\frac12 Q\T mQ+ M \la{barM} \ee 
is a momentum-like quantity. 
\la{soptex}
\end{corollary}
\noindent {\em Proof}: The proof of this statement can be carried out in two 
stages. In the first stage, we consider an extremal solution pair $(Q,P)\in W^k_m
\subset W^k_{n,N}$ and obtain an expression for $P$ in terms of $Q$ and $M$. 
We observe that if 
\[ P= \mu Q+QM, \]
where $\mu\in\bbr^{nn}$ is constant and $M$ is the body momentum, then 
$(Q,P)\in W^k_{n,N}$. This can be easily verified by taking a time derivative of 
$P$ along the extremal solutions given by equations (\ref{sOCP}) and (\ref{Madot}). 
We also observe that $PQ\T=\mu+ m$ for this solution pair. This solution is in 
$W^k_m\subset W^k_{n,N}$ if and only if 
\beas PQ\T+QP\T= \mu+\mu\T= k, \;\ &\mbox{and }& \\
PQ\T-QP\T= \mu-\mu\T+2m =m \ &\Leftrightarrow &\ \mu\T-\mu= m. \eeas 
This gives $\mu$ uniquely as 
\[ \mu =\frac12 (k-m), \]
and thus, $(Q,P)\in W^k_m$ where 
\be P= \frac12 (k-m)Q+ QM. \la{QPinWkm} \ee
Now applying Proposition \ref{Wkto0} in the second stage of this proof, we have 
$(Q,P_0)\in W^0_m$ where $P_0=\bar{P}-\frac12 kQ$, which when applied to 
(\ref{QPinWkm}) gives the result expressed in equation (\ref{Peqn}). As a result of 
Proposition \ref{Wkto0}, we also know that the momenta $m$ and $M$ are left 
unchanged by this transformation in the costate variable. That $P_0=Q\bar{M}$ is 
easily verified by susbtituting $\bar{M}$ from equation (\ref{barM}) and comparing 
with equation (\ref{Peqn}). \qed \\
Note that, with the costate variable $P=P_0$, the solution to the optimal control 
problem satisfies
\be k= QP_0\T+ P_0Q\T= -\frac12 m+ \frac12 m= 0. \la{QPrel} \ee
Also note that, for the special case of the rigid body in $N$ 
dimensions ($n=N$), the momentum quantity $\bar{M}=\frac12 M$ is half the body 
momentum. \\
%The $(Q,P_0)$ variables may 
%be considered as a symmetric representation of this Hamiltonian system in $V(n,N)\times 
%V(n,N)$, where we relate $W^0_{n,N}$ to $V(n,N)\times V(n,N)$ by
%\[ (P_0 Q\T -s)Q= c\T Q\in V(n,N), \] where $s$ given by (\ref{sP}) and 
%$c\in SO(n)$ given by (\ref{pBeq}) are constant along the flow of this system. \\ 
 
This result has further important implications for the symplectic structure of 
$W^0_{n,N}$ and its relation to the symplectic structure on  
$T^\star V(n,N)$, which we will explore in the next section.
Note that the equations (\ref{sOCP})-(\ref{solOCP}) conserve $PQ\T$ and $QP\T$ 
separately. In the symmetric representation of the rigid body equations given in 
Bloch \etal (2002), $PQ\T$ and hence $PQ\T- QP\T$ and $PQ\T+ QP\T$ are constant, 
and this generalizes to the extremal solution (in $W^k_{n,N}$) of the optimal 
control problem on the Stiefel manifold $V(n,N)$ for $1\le n<N$. 

\subsection{Correspondence of the Variational and Optimal Control Solutions}

We now give the correspondence between the variational (or Lagrangian) and optimal 
control (or Hamiltonian) respresentations of the extremal solutions to this problem. 
The extremal solutions to the variational problem can be defined in terms of the 
pair $(Q,\dot Q)\in TV(n,N)$ or $(Q,S)\in T^\star V(n,N)$ and extremal solutions to 
the optimal control problem are defined in terms of the pair $(Q,P_0)\in W^0_{n,N}$, 
a symplectic manifold. Here we describe a correspondence between them, $\mfrak{M}: 
T^\star V(n,N)\rightarrow W^0_{n,N}$ such that $\mfrak{M}:\ (Q,S) \mapsto (Q,P_0)
\in W^0_{n,N}$. \\ 
%\footnote{pronounced: m\`{a}} {\bngxii m}: $T^\star V(n,N)\rightarrow W^0_{n,N}$ such 
%that {\bngxii m} $:\ (Q,S) \mapsto (Q,P_0)\in W^0_{n,N}$. \\ 

From the solutions to the variational and optimal control problems, we know that the 
body momentum $M$ satisfies
\be M= Q\T S- S\T Q= Q\T P_0- P_0\T Q. \la{khokhla} \ee
In addition, on pre-multiplying equation (\ref{khokhla}) by $Q$ and post-multiplying 
it by $Q\T$ on both sides, we also get 
\be m= QMQ\T= SQ\T- QS\T= P_0Q\T- QP_0\T, \la{smkhok} \ee 
which is derived from (\ref{khokhla}). Equation (\ref{khokhla}) leads us to the 
following result. 
\begin{corollary}
The extremal solutions $(Q,S)$ and ($Q,P_0)$ to the variational and optimal control 
problems, respectively, are related by 
\be P_0= S+X \la{PSrel} \ee
where 
\be X= DQ,\;\ D=D\T. \la{Xform} \ee
\la{P02S}
\end{corollary} 
\noindent This statement is a corollary of Proposition \ref{LQker}, and is obtained 
from inspection of equation (\ref{khokhla}). The difference $P_0-S=X$ in equation 
(\ref{PSrel}) obviously lies in the kernel of the linear map $L_Q : \bbr^{nN} 
\rightarrow \bbr^{NN}$ defined in (\ref{LQmap}). From Proposition \ref{LQker}, we 
know that $X$ is going to have the form given by equation (\ref{Xform}). \\

For the extremal solution to the optimal control problem $(Q,P_0)\in W^0_m\subset
W^0_{n,N}$, we know from Corollary \ref{soptex} that $P_0=-\frac12 mQ+QM$. The 
extremal solution of the variational problem in Section \ref{sec:varp} was given in 
terms of $(Q,S)\in T^\star V(n,N)$. The following result gives the relation between 
these two solution pairs. 
\begin{proposition}
The map %{\bngxii m} 
$\mfrak{M} :T^\star V(n,N)\rightarrow W^0_{n,N}$ defined by %{\bngxii m}
$\mfrak{M}((Q,S))=(Q,P_0)$, is a diffeomorphism and $P_0$ is given in terms of $Q$ 
and $S$ by  
\be P_0= S+DQ, \la{PSrel1} \ee 
where  
\be D=D\T= \frac12 Q[\Lam,U] Q\T= -\frac12 (QS\T+SQ\T). \la{Dform} \ee 
The inverse of this diffeomorphism is given by %{\bngxii m}
$\mfrak{M}^{-1}: (Q,P_0)\mapsto (Q,S)$ where 
\be S= P_0-[(Q\Lam^{-1}P_0\T)(Q\Lam^{-1}Q\T)^{-1}-Q\Lam^{-1}Q\T\wh{J_Q}^{-1}(R_0)]Q, 
\la{SfP0} \ee
where $\wh{J_Q}^{-1}(R_0)$ is given by equations (\ref{Req}) with $S$ replaced by 
$P_0$. 
%\bea U_1&=& \Lam^{-1}\big(Q\T(Q\Lam^{-1}Q\T)^{-1}P_0 -P_0\T (Q\Lam^{-1}Q\T)^{-1}Q\big)
%\Lam^{-1}, \nn \\ 
%U_2 &=& \Lam^{-1}(Q\T \wh{J_Q}^{-1}(R) Q)\Lam^{-1}, \la{U1U2} \eea
%and 
%\bea R&=& (Q\Lam^{-1}P_0\T)(Q\Lam^{-1}Q\T)^{-1}-(Q\Lam^{-1}Q\T)^{-1}(P_0\Lam^{-1}Q\T), 
%\nn \\ \wh{J_Q}(X)&=& (Q\Lam^{-1}Q\T) X+ X(Q\Lam^{-1} Q\T);\;\ \wh{J_Q}: \mfrak{so}(n)
%\rightarrow \mfrak{so}(n) \la{Req} \eea
\la{vroprel}
\end{proposition} 
\noindent {\em Proof}: From Corollary \ref{P02S}, we know that $P_0= S+DQ$ where $D$ is 
symmetric. Thus, we have 
\[ S=P_0-DQ= -(\frac12 m+ D)Q+ QM, \]
where $M=Q\T S-S\T Q$ and $m=QMQ\T$. This gives us 
\beas -(\frac12 m+ D)Q+ Q(Q\T S-S\T Q)&=& S \\
%\Leftrightarrow -(\frac12 m+ D)Q+ S-QS\T Q&=& S \\
\Leftrightarrow -\frac12 (SQ\T-QS\T)Q- DQ&=& QS\T Q \\ 
%\Leftrightarrow -\frac12 SQ\T Q- DQ&=& \frac12 QS\T Q \\
\Leftrightarrow - DQ&=& \frac12 (QS\T+SQ\T) Q \\
\Leftrightarrow D= -\frac12 (QS\T+SQ\T)&=& -\frac12 Q[\Lam,U] Q\T \eeas
since $S= QU\Lam$, where $QU$ is obtained from $Q$ and $S$ using Theorem \ref{algeq}. 
This establishes the relation given by equations (\ref{PSrel1})-(\ref{Dform}) 
between the variables $S$ and $P_0$, and $P_0$ is seen to be linearly dependent on 
$S$. Clearly, if $S=0$ then $P_0=0$. For the converse, if $P_0=0$ then $S=-DQ$. Since  
$D$ is symmetric, we have $Q\T S=S\T Q= -Q\T DQ$. Since $S=QU\Lam$, we have $Q\T S-
S\T Q= J_Q(U)$, as defined in equation (\ref{JQU}). From Lemma \ref{JQlem}, we know 
that if $J_Q (U)=0$, then $QU=0$. Thus, when $P_0=0$, we have $J_Q(U)=0$ and hence 
$S=QU\Lam=0$. Thus %{\bngxii m}$
$\mfrak{M}((Q,S))=(Q,P_0)$ as given by (\ref{PSrel1}) is a 
diffeomorphism. For the inverse of this map, it is clear from Theorem \ref{algeq} and 
its proof that $U$ can be expressed in terms of $Q$ and $P_0$ in the same manner that 
it is expressed in terms of $Q$ and $S$ in that theorem. This is true because that 
representation is based on the decomposition $M=Q\T S-S\T Q$, which has the same form 
as $M=Q\T P_0-P_0\T Q$. Hence, we also have a diffeomorphism $\mcal{Z}_0: W^0_{n,N}
\rightarrow T V(n,N)$ defined by $\mcal{Z}_0((Q,P_0))=(Q,[U])=(Q,[U_1+U_2])$ as follows
\beas U_1 &=& \Lam^{-1}\big(Q\T(Q\Lam^{-1}Q\T)^{-1}P_0 -P_0\T (Q\Lam^{-1}Q\T)^{-1}Q\big)
\Lam^{-1}, \nn \\ U_2 &=& \Lam^{-1}(Q\T \wh{J_Q}^{-1}(R_0) Q)\Lam^{-1}, \eeas 
where 
\[ R_0=(Q\Lam^{-1}P_0\T)(Q\Lam^{-1}Q\T)^{-1}-(Q\Lam^{-1}Q\T)^{-1}(P_0 \Lam^{-1}Q\T). \] 
Thus the inverse of %{\bngxii m} 
$\mfrak{M}$ is given by %{\bngxii m}$
$\mfrak{M}^{-1}= \mcal{Z}^{-1}\circ\mcal{Z}_{0}$, where $\mcal{Z}$ is as defined 
in Corollary \ref{S2Uc}. From this corollary, we also know that $S=Q(U_1+U_2)\Lam$ 
uniquely determines $S$, and this gives us the relation (\ref{SfP0}) for $S$ in terms 
of $P_0$ and $Q$. \qed \\

Note that if $\Lam=I_N$, the $N\times N$ identity matrix, then $S=P_0$. Also note 
that $D=(Q\Lam^{-1}P_0\T)(Q\Lam^{-1}Q\T)^{-1}-Q\Lam^{-1}Q\T\wh{J_Q}^{-1}(R_0)$ in 
(\ref{SfP0}) is a symmetric matrix. This can be shown using equations (\ref{Req}) of 
Theorem \ref{algeq}, as follows
\beas \fl D- D\T = (Q\Lam^{-1}P_0\T)(Q\Lam^{-1}Q\T)^{-1}-Q\Lam^{-1}Q\T\wh{J_Q}^{-1}
(R_0)-\wh{J_Q}^{-1}(R_0)Q\Lam^{-1}Q\T \\  -(Q\Lam^{-1}Q\T)^{-1}(Q\Lam^{-1}P_0\T) \\ 
\fl = (Q\Lam^{-1}P_0\T)(Q\Lam^{-1}Q\T)^{-1}-Q\Lam^{-1}Q\T X_0- X_0 Q\Lam^{-1}Q\T-
(Q\Lam^{-1}Q\T)^{-1}(Q\Lam^{-1}P_0\T) \\ 
\fl = R_0- \wh{J_Q}(X_0)= R_0- R_0=0, \eeas
where $X_0= \wh{J_Q}^{-1}(R_0)$. This expresses the symmetric matrix $D$ in terms of 
$Q$ and $P_0$, while $D$ is also given as a function of $Q$ and $S$ by equation 
(\ref{Dform}).

\section{The tangent and cotangent bundles of the Stiefel Manifold} 

In this section, we explore the structure of the tangent and cotangent bundles of the 
Stiefel manifold, which is the homogeneous space $V(n,N)=SO(N)/SO(N-n)$. The Stiefel 
manifold is parametrized by $Q\in\bbr^{nN}$ such that $QQ\T=I_n$. Consider the point 
$Q_0=[I_n\;\ 0]$ on $V(n,N)$, where $0$ denotes the $(N-n)\times n$ matrix of zeros. 
$SO(N)$ acts on $Q_0$ on the right as $\mcal{R}:Q_0\mapsto Q_0R$ where $R\in SO(N)$, 
and the isotropy group is the subgroup 
\[ \Big\{R\in SO(N)\ |\ R=\matl{cc} I_n & 0\\ 0& \bar{R} \matr;\ \bar{R}\in SO(N-n) 
\Big\}. \]
This gives the dimension of $V(n,N)$ as the difference of the dimensions of $SO(N)$ and 
$SO(N-n)$
\[ \mbox{dim}\ V(n,N)= N(N-1)/2- (N-n)(N-n-1)/2= nN-n(n+1)/2. \]
%as we obtained earlier using the constraint $QQ\T=I_n$. \\

\subsection{Symplectic Structure of the cotangent bundle of the Stiefel Manifold} 

The tangent space to $V(n,N)$ at $Q$ is parametrized by $QU$, $U\in\mfrak{so}(N)$. 
However, $U$ is unique only upto the equivalence class $[U]$ defined by equation 
(\ref{eqvclass}). If $Q=[I_n\;\ 0]$ as before, then the set $V\in SO(N)$ such that 
$QV=0$ is given by 
\[ \left\{ V\in SO(N)\ \big|\ V=\matl{cc} 0 & 0\\ 0 & \bar{V} \matr,\;\ \bar{V}\in 
\mfrak{so}(N-n)\right\}. \]
This determines $T_Q V(n,N)$ as the vector space 
\[ \mfrak{so}(N)/\mfrak{so}(N-n). \]
We occasionally refer to elements of $T_Q V(n,N)$ as the equivalence classes $Q[U]$, 
where $[U]$ is as defined by equation (\ref{eqvclass}). \\

We may parametrize the co-tangent space $T^{\star}_Q V(n,N)$ by 
\[ P\in T^{\star}_Q V(n,N),\;\ P=2\langle QM, \cdot\rangle,\;\ M\in \mfrak{so}(N). \] 
Let $X\in T_Q V(n,N)$, hence $X=QU$, $U\in\mfrak{so}(N)$. Then we have 
\[ P(X)= 2\langle QM, QU\rangle= 2\langle Q\T QM, U\rangle. \]
In the case $n=N$, this pairing gives the Killing form on $\mfrak{so}(N)$, which 
is non-degenerate. The following result shows that in the equivalence classes 
$[\cdot]$ defined on $\mfrak{so}(N)$, this pairing between $T_Q V(n,N)$ and 
$T^{\star}_Q V(n,N)$ is non-degenerate.
\begin{lemma}  
The pairing %\footnote{pronounced: \^{o}\`{i}} {\bngxii OI} $: 
$\mfrak{O}: T^{\star}_Q V(n,N)\times T_Q V(n,N) \rightarrow\bbr$ given by %{\bngxii OI}$
$\mfrak{O}(Q[M],Q[U])\mapsto \langle Q[M],Q[U]\rangle$ is non-degenerate. 
\la{nondp}
\end{lemma}
\noindent{\em Proof}: For this pairing to be non-degenerate on $T_Q V(n,N)$, we should 
have \[ \langle QM,QU\rangle= 0\ \forall\ U\ \Rightarrow\ Q[M]=0. \]
Evaluating the pairing on the left hand side gives us 
\beas \langle QM,QU\rangle= \langle Q\T QM,U\rangle&=&0\ \forall\ U\in\mfrak{so}(N) \\
\Leftrightarrow Q\T QM+ MQ\T Q &=& 0. \eeas
Now, by Lemma \ref{JQlem}, the above equality is satisfied if and only if $QM=0$ 
(substituting $U=M$ and $\Lam= I_N$ in Lemma \ref{JQlem}). 
%Alternate proof :-
%To show this, we again consider $Q_0=[I_n\;\ 0]\in V(n,N)$, 
%since we can obtain any other point in $V(n,N)$ by the left action of $SO(n)$, 
%and the pairing is invariant to this left action. In this case, if 
%\[ M=\matl{cc} M_1 & M_2\\ -M_2\T & M_3 \matr,\;\ U=\matl{cc} U_1 & U_2\\ -U_2\T & U_3 
%\matr, \]
%where $M_1, U_1\in \mfrak{so}(n)$, $M_3, U_3\in\mfrak{so}(N-n)$, and $M_2, U_2\in 
%\bbr^{n(N-n)}$, we have $Q_0M=[M_1\;\ M_2]$ and $Q_0U=[U_1\;\ U_2]$. Thus, the above 
%pairing takes the form 
%\beas \langle Q_0M, Q_0U\rangle&=& \Tr (M_1\T U_1)+\Tr (M_2\T U_2) \\
%&=& \langle M_1,U_1\rangle+ \langle M_2, U_2\rangle. \eeas 
%For any element of $[M]$, the $M_1$ and $M_2$ matrices are identical, and for any 
%element of $[U]$, the $U_1$ and $U_2$ matrices are identical. Since the above pairing 
%\[ \big((M_1,U_1),(M_2,U_2)\big)= \langle M_1,U_1\rangle+ \langle M_2,U_2\rangle \]
%is non-degenerate, 
Thus, we have shown that the bilinear pairing %{\bngxii OI}$
$\mfrak{O}(Q[M],Q[U])\mapsto \langle Q[M],Q[U]\rangle$ is 
non-degenerate on $T_Q V(n,N)$. \qed \\ 
Hence, this representation of $T^{\star}_Q V(n,N)$ is well defined. \\

Identifying $T^\star V(n,N)$ and $T V(n,N)$ by this pairing, we may parametrize 
$TT^\star V(n,N)$ by vectors $(Q,QM,QU,QUM+QZ)$, where $M, U, Z\in \mfrak{so}(N)$. 
Viewing the co-tangent bundle $T^\star V(n,N)$ as a subset of $\bbr^{nN}\times
\bbr^{nN}$, we may pull back the symplectic form on $\bbr^{nN}\times\bbr^{nN}$ given by 
\[ \omega((A_1,B_1),(A_2,B_2))= \langle B_2, A_1\rangle -\langle B_1,A_2\rangle, \]
to $T^\star V(n,N)$ via the inclusion map. If 
\[ \fl X_1= (Q,QM,QU_1,QU_1 M+QZ_1),\;\ X_2= (Q,QM,QU_2,QU_2 M+QZ_2)\in T_{(Q,QM)}T^\star 
V(n,N), \]
then the two form we obtain on $T^\star V(n,N)$ is 
\bea \fl \omega_{(Q,QM)}(X_1,X_2)&=& \langle QU_2 M+QZ_2, QU_1\rangle- \langle QU_1 M
+QZ_1, QU_2\rangle \nn \\ \fl &=& \langle Q\T QU_1,Z_2\rangle- \langle Q\T QU_2,Z_1
\rangle+ \langle M, U_1Q\T QU_2-U_2 Q\T QU_1 \rangle. \la{symfpb} \eea 
It is simple to check using Lemma \ref{JQlem} that this is indeed nondegenerate and 
hence a symplectic form. This expression for the symplectic form can also be obtained 
using the canonical structure on $T^\star V(n,N)$. In general, for a smooth manifold 
$M$ and using the projection $\pi: T^\star M\rightarrow M$, we define a one form 
$\theta$ on $T^\star M$ by 
\[ \theta_p (X)= p\pi_{\star} X, \]
where $p\in T^\star M$ and $X$ is a vector field on $T^\star M$. We define the canonical 
symplectic form $\omega_c$ on $T^\star M$ by setting $\omega_c= -d\theta$, where 
$d\theta$ is the exterior derivative of the one form $\theta$. We may simplify this 
expression using the identity (given in Bloch \etal (2003))
\[ d\theta (X,Y)= X(\theta(Y))-Y(\theta (X))-\theta ([X,Y]), \]
so 
\be \omega_c (X,Y)= Y(\theta (X))-X(\theta (Y))+\theta ([X,Y]). \la{canform} \ee

We apply equation (\ref{canform}) to the case where the manifold $M=V(n,N)$. 
We parametrize $T^\star V(n,N)$ by pairs $(Q,M)$, $Q\in V(n,N)$ and $M\in 
\mfrak{so}(N)$. Then we can write elements of $T_{(Q,M)} T^\star V(n,N)$ 
at the point $(Q,M)$ as 
\[ X_k= (Q U_k, Z_k),\;\ k=1,2. \] 
Hence, if $\pi$ is the projection $\pi: TT^\star V(n,N) \rightarrow T V(n,N)$, 
then $\pi_\star X_k= QU_k$, and if $p\in T^\star V(n,N)$, $p= \langle QM,\cdot 
\rangle$, then 
\[ \theta (X_k)= p\pi_\star X_k= \langle QM, QU_k \rangle. \] 
Thus, from equation (\ref{canform}) we have
\be \omega_c (X_1,X_2)= X_2\langle QM,QU_1\rangle- X_1\langle QM,QU_2\rangle
+ \langle QM,\pi_\star [X_1,X_2]\rangle. \la{omc} \ee
The last term in the above equation can be simplified as follows
\beas \pi_\star [X_1, X_2]&=& [\pi_\star X_1, \pi_\star X_2]\circ\pi \\
&=& [QU_1, QU_2] \;\ \mbox{(as vector fields)} \\
&=& QU_1 U_2-QU_2 U_1 \\
&=& Q[U_1, U_2]. \eeas
The first two terms are given by
\[ \fl X_2\langle QM,QU_1\rangle = (QU_2,Z_2)\langle QM,QU_1\rangle 
= \langle QZ_2, QU_1\rangle+ \langle QU_2 M,QU_1\rangle+\langle QM, QU_2 
U_1\rangle, \] and 
\[ \fl X_1\langle QM,QU_2\rangle = (QU_1,Z_1)\langle QM,QU_2\rangle 
= \langle QZ_1, QU_2\rangle+ \langle QU_1 M,QU_2\rangle+\langle QM, QU_1 
U_2\rangle. \] 
Thus, we obtain
\bea \omega_c(X_1,X_2)&=& \langle QZ_2, QU_1\rangle- \langle QZ_1, QU_2\rangle
-\langle M,U_2Q\T QU_1\rangle+ \langle M,U_1Q\T QU_2\rangle \nn \\ & &+\langle M,
Q\T Q(U_2 U_1-U_1 U_2)\rangle +\langle QM,Q[U_1,U_2]\rangle \nn \\ &=& \langle QZ_2, 
QU_1\rangle- \langle QZ_1, QU_2\rangle+\langle M,U_1Q\T QU_2-U_2Q\T QU_1\rangle, 
\la{omccanf} \eea
which is identical to equation (\ref{symfpb}) obtained from restricting the 
two-form on $\bbr^{nN}\times\bbr^{nN}$ to $T^\star V(n,N)$. \\

We now have a formula for the natural symplectic form on $T^\star V(n,N)$, and 
know that it can be formulated as the restriction of the symplectic form on the 
product $\bbr^{nN}\times\bbr^{nN}$. We want to use this formula to recover 
the Hamiltonian flow corresponding to the geodesic problem (\ref{OCPnN}) on 
$V(n,N)$. The Hamiltonian is given by (\ref{Ham}) as $H(P,Q,U)= \langle P,QU \rangle
-\langle QU, QU\Lam\rangle$ and the optimal control $U^*$ is given by the maximum 
principle as in (\ref{maxprin}). We replace $P$ by the parametrization $QM$, where 
we expect $M$ to be the momentum. Thus, we can write the Hamiltonian as 
\be H(M,Q,U)= \langle QM,QU\rangle-\frac12\langle QU\Lam,QU\rangle. 
\la{MHam} \ee
The Hamiltonian flow on $T^\star V(n,N)$ is obtained from the solution of 
\[ dH_{(M,Q,U)} X_2= \omega_c (X_1,X_2),\;\ \forall\ X_2=(QU_2,Z_2), \]
where $X_1=X_H$ is the Hamiltonian vector field corresponding to $H$. We calculate 
\[ \fl dH (X_2)= \langle QU_2 M,QU\rangle+ \langle QM,QU_2 U\rangle +\langle QZ_2, 
QU\rangle -\frac12\langle QU_2 U\Lam,QU\rangle -\frac12\langle QU\Lam,QU_2 U\rangle, \]
and this is equated to $\omega_c (X_H,X_2)$ in (\ref{omccanf}). From equating these 
expressions after replacing $Z_1$ and $U_1$ by $Z$ and $U$, respectively, in 
(\ref{omccanf}), we get
\be \langle QU_2,QZ\rangle= \langle QU_2,QMU \rangle- \langle QU_2,QU M
\rangle-\langle QU_2,QU\Lam U\rangle. \la{poch} \ee
Note that $\langle QU_2 Q\T,QU\Lam UQ\T\rangle=0$ since $QU_2 Q\T$ is skew-symmetric 
and $QU\Lam UQ\T$ is symmetric. Therefore, $\langle QU_2, QU\Lam UQ\T Q\rangle=0$ 
and we can express the last term in equation (\ref{poch}) as  
\[ \langle QU_2,QU\Lam U\rangle=\langle QU_2,Q(Q\T QU\Lam U-U\Lam UQ\T Q)\rangle. \]
If $F$ is skew-symmetric and $\langle QU_2,QU\Lam U\rangle=\langle QU_2,QF\rangle$, 
then 
\[ \langle QU_2,Q(F-[Q\T Q,U\Lam U])\rangle= \langle QU_2,QU\Lam UQ\T Q\rangle=0\;\ 
\forall\ U_2. \]
Thus $F=[Q\T Q,U\Lam U]$ up to an equivalence class. % defined by the above inner 
%product relation, and $[Q\T Q,U\Lam U]$ is the unique skew-symmetric 
%extension of $QU\Lam U$ in equation (\ref{poch}). 
Hence, from (\ref{poch}) we get 
\beas \lefteqn{\langle QU_2,Q(Z_1-[M,U]+[Q\T Q, U\Lam U])\rangle=0\ \forall\ U_2} \\
\lefteqn{\Rightarrow Z_1= [M,U]-[Q\T Q,U\Lam U],} \eeas
up to equivalence class. It follows that the geodesics on $T^\star V(n,N)$ with the 
metric in (\ref{OCPnN}) are given by the Hamiltonian flow
\beas \dot Q&=& QU,\\
\dot M&=& [M,U]-A,\;\ A=[Q\T Q,U\Lam U], \eeas
as given previously by (\ref{VPsol}) and (\ref{Madot}).  

\subsection{Symplectomorphism between the cotangent bundle and $W^0_{n,N}$}

Finally we explore the relation between the symplectic structures on 
the cotangent bundle $T^\star V(n,N)$ and $W^0_{n,N}$. %The symplectic two-form 
%is evaluated on $W^0_{n,N}$ at the extremal solution pair $(Q,P_0)\in W^0_{n,N}$, 
Consider the Hamiltonian vector fields defined by (\ref{sOCP})-(\ref{solOCP}):
\beas X_1= (QU_1, P_0U_1-QA_1), \\
X_2= (QU_2, P_0U_2-QA_2), \eeas 
where $A_i= [Q\T Q, U_i\Lam U_i]$, $i=1,2$. The symplectic form $\Omega$ on 
$W^0_{n,N}$ evaluated along these vector fields on $W^0_{n,N}$ is obtained using 
equation (\ref{2fWnN}) as: 
\be \Omega_{(Q,P_0)} (X_1,X_2) = %&=& \langle \partial_2 P, \partial_1 Q \rangle
%- \langle \partial_1 P, \partial_2 Q \rangle \nn \\
\langle P_0U_2- QA_2,QU_1 \rangle- \langle P_0U_1-QA_1, QU_2 \rangle. \la{VnNVnN} \ee 
From Corollary \ref{soptex}, we know that if $(Q,P_0)\in W^0_m\subset W^0_{n,N}$, 
then 
\[ P_0=-\frac{m}{2}Q+QM = Q\bar{M}, \]
where $\bar{M}$ is given by (\ref{barM}). 

Now we show the following relationship between $T^\star V(n,N)$ and $W^0_{n,N}$. 
\begin{theorem}
The map $\Phi: (\Omega, W^0_{n,N}) \rightarrow (\omega_c, T^\star V(n,N))$, given 
by $(Q,P_0) \mapsto (Q,Q\bar{M})$ where 
\be \bar{M}=-\frac12 Q\T(P_0 Q\T- QP_0\T)Q+ Q\T P_0-P_0\T Q, \la{bMP0} \ee 
is a symplectomorphism. 
\la{W0Ts}
\end{theorem} 
\noindent {\em Proof:} Coordinates for $T^\star_Q V(n,N)$ are given by 
$(Q,Q\bar{M})$. Consider the Hamiltonian vector fields on $T^\star V(n,N)$ 
corresponding to the extremal flows for the optimal control problem considered 
in the last section
\beas \chi_1=\partial_1 (Q,Q\bar{M})= (QU_1, QU_1 \bar{M}+QZ_1), \\
\chi_2= \partial_2 (Q,Q\bar{M})= (QU_2, QU_2 \bar{M}+QZ_2), \eeas 
where 
\beas Z_1= \partial_1 \bar{M}= [\bar{M},U_1]-A_1, \\
Z_2= \partial_2 \bar{M}= [\bar{M},U_2]-A_2. \eeas 
We show that the push-forwards of the Hamiltonian vector fields $X_i$ on $W^0_{n,N}$ 
defined earlier, give the Hamiltonian vector fields $\chi_i$ on $T^\star V(n,N)$. 
%are given by $\Phi_\star X_i (\Phi(Q,P_0))= \chi_i (Q,Q\bar{M})$, 
%since $\chi_i (Q,Q\bar{M})= (QU_i, Q\bar{M}U_i-QA_i)$ and $Q\bar{M}=P_0$. 
We have 
\beas \big(\Phi_\star X_i\big)(\Phi(Q,P_0))&=& \big(T_{(Q,P_0)}\Phi\cdot X_1\big)
(\Phi(Q,P_0))\\ &=& (QU_i, Q\bar{M}U_i-QA_i) \\ &=& (QU_i, QU_i\bar{M}+Q[\bar{M},
U_i]-QA_i) \\ &=& (QU_i, QU_i\bar{M}+QZ_i)= \chi_i (Q,Q\bar{M}) \\ 
\Leftrightarrow \Phi_\star X_i &=& \chi_i. \eeas
The canonical symplectic form on $T^\star V(n,N)$ is evaluated as 
\bea \fl \omega_{c\ (Q,Q\bar{M})} (\chi_1,\chi_2) &=& \langle QU_2 \bar{M}+QZ_2, 
QU_1 \rangle- \langle QU_1 \bar{M}+QZ_1, QU_2 \rangle \nn \\ \fl 
&=& \langle QU_2 \bar{M}+ Q[\bar{M},U_2]-QA_2, QU_1 \rangle- \langle
QU_1 \bar{M}+ Q[\bar{M},U_1]-QA_1, QU_2 \rangle \nn \\ \fl
&=& \langle Q\bar{M} U_2- QA_2, QU_1 \rangle- \langle Q\bar{M} U_1
- QA_1, QU_2 \rangle. \la{TsVnN} \eea 
The pull-back of this symplectic form from $T^\star V(n,N)$ to $W^0_{n,N}$ 
gives us the symplectic form on $W^0_{n,N}$, as shown below:
\beas \fl \left(\Phi^\star \omega_c\right)_{(Q,P_0)} (X_1,X_2)&=& 
\omega_{c\ (Q,Q\bar{M})} (\Phi_\star X_1,\Phi_\star X_2)= 
\omega_{c\ (Q,Q\bar{M})} (\chi_1,\chi_2) \\ \fl &=& \langle Q\bar{M}U_2- QA_2, 
QU_1 \rangle- \langle Q\bar{M} U_1- QA_1, QU_2 \rangle \\ \fl
&=& \langle P_0 U_2-QA_2, QU_1 \rangle- \langle P_0 U_1-QA_1,QU_2 \rangle = 
\Omega_{(Q,P_0)}(X_1,X_2). \eeas 
However, from Lemma \ref{LiWknN}, we know that the vector fields $X_i$ form 
a Lie algebra that spans the tangent space to $W^0_{n,N}$ at every point $(Q,P_0)$. 
Hence, we have
\be \Phi^\star\omega_c= \Omega. \la{PbTs2SVnN} \ee
The map $\Phi: W^0_{n,N} \rightarrow T^\star V(n,N)$ 
given by $(Q,P_0) \mapsto (Q,Q\bar{M})$ where $\bar{M}$ is given by (\ref{bMP0}),  
%\[ \bar{M}= -\frac12 Q\T(P_0 Q\T- QP_0\T)Q+ Q\T P_0-P_0\T Q, \]
has an inverse which is simply given by $\Phi^{-1} :(Q,Q\bar{M}) \mapsto 
(Q,P_0)$ where $P_0=Q\bar{M}$. The inverse map $\Phi^{-1}$ is clearly a 
diffeomorphism. Thus, the map $\Phi$ is a diffeomorphism that also maps the 
symplectic form in $W^0_{n,N}$ to the symplectic form in $T^\star V(n,N)$, and is 
hence a symplectomorphism. \qed \\ 

With $\bar{M}$ expressed in terms of $Q$ and $P_0$ as in this theorem, we have  
that $Q\bar{M}=P_0$ where $(Q,P_0)\in W^0_{n,N}$. Hence, Theorem \ref{W0Ts} shows 
that the solution space $W^0_{n,N}$ for the optimal control problem can be 
identified with the cotangent bundle of the Stiefel manifold. 
The Hamiltonian formulation of the optimal control problem parametrized by 
$(Q,P_0)$ or $(Q,Q\bar{M})$ is related to the variational formulation of this 
problem parametrized by $(Q,\dot Q)$ or $(Q,S)$, via the relationship established 
in Proposition \ref{vroprel}. Combining Theorem \ref{W0Ts} with Proposition 
\ref{vroprel}, we get the following expression for $\bar{M}$ in terms of $(Q,S)$:
\be \bar{M}= Q\T\bar{S}-\bar{S}\T Q,\;\ \bar{S}=S(I_N-\frac12 Q\T Q). \la{QS2bM} \ee
Of course, in the case case $\Lam=I_N$, we have $P_0=S$ and these solutions are 
directly related. 
%\begin{corollary} 
%The fiber transformation $\mathbb{F}\mcal{L}=$ {\bngxii m} $\circ\ \Phi: TV(n,N)
%\rightarrow T^\star V(n,N)$ for the Lagrangian $\mcal{L}$ defined in (\ref{Lag}) 
%is the Legendre transform, $\mathbb{F}\mcal{L} (Q,S)=(Q,Q\bar{M})$. 
%\end{corollary} 
%\noindent The Legendre transform $\mathbb{F}\mcal{L}$ relates the extremal solution 
%pairs $(Q,S)\in TV(n,N)$ and $(Q,Q\bar{M})\in T^\star V(n,N)$ which were obtained 
%from the variational and optimal control formulations, respectively.

\section{Applications and Open Problems}

Nowadays numerical linear algebra computations and numerical integration of 
ODEs are increasingly based on variational problems on manifolds. Problems on 
Stiefel (and Grassman) manifolds are finding increasing use in numerical 
linear algebra applications (see, for example, Edelman \etal (1998), Elden and 
Park (1999), and references therein). The first of these papers develops Newton 
and conjugate gradient methods on these manifolds, while the latter paper deals 
with a problem related to regression analysis in psychometrics. Another application 
of numerical calculations on Stiefel manifolds is in computing Lyapunov exponents 
for finite-dimensional dynamical systems by time integration. The Lyapunov exponents 
are computed by a continuous orthonormalization (which is essential for stable 
numerical integration) of a set of solution vectors of the linearized system. This 
amounts to restricting the linearized system to the Stiefel manifold $V(k,N)$ for 
computing the $k$ largest Lyapunov exponents of an $N$-dimensional system. A 
sample of the literature on this topic can be obtained from Bridges and Reich (2001) 
and references therein. \\

It is known that the geodesic flow on the Stiefel manifold $V(n,N)$ with a 
left-invariant metric is integrable for the extreme cases: when $n=1$, which 
represents the geodesic flow on the sphere/ellipsoid; and when $n=N$, which 
represents the $N$-dimensional rigid body on $SO(N)$. Bolsinov and Jovanovic have 
shown that the extremal flows on Stiefel manifolds and other homogeneous spaces 
with bi-invariant metrics are integrable. There are two sets of integrals for 
such flows. The first set of integrals are the Noether integrals 
\[ \mathcal{F}_1=\{ h\circ\mfrak{G},\ %\mbox{\bngxii NG},\ \mbox{\bngxii NG}: 
\mfrak{G}: TV(n,N)\rightarrow \mfrak{so}^\star(N),\ h: \mfrak{so}^\star(N)
\rightarrow \bbr\},\] 
where %{\bngxii NG}$
$\mfrak{G}(Q,QU)=J_Q (U)=M$. In this case, it is easy to verify that 
$M$ is conserved along the extremal flows. Without loss of generality, we may 
take $\Lam=I_N$ as the bi-invariant metric, and the extremal flows are then 
given by 
\[ \dot Q=QU,\;\ \ddot Q= BQ, \ B=B\T. \] 
The body momentum is $M=Q\T QU+UQ\T Q$, and its derivative along the extremal 
flows is 
\[ \dot M= \dot Q\T\dot Q+Q\T\ddot Q-\ddot Q\T Q-\dot Q\T\dot Q=0. \] 
This is a generalization of the case of the symmetric rigid body, where $\Lam$ 
is a scalar multiple of the identity matrix. The second set of integrals for 
this bi-invariant case is  
\[ \mathcal{F}_2=\mbox{all $SO(N)$-invariant functions on $TV(n,N)$}. \]

The most important issue for future research is the integrability of the 
extremal solutions in the general case $(1< n< N)$ with a left-invariant metric. 
Another research issue of considerable interest for numerical applications 
is discretization of the optimal control problem (\ref{OCPnN}) based on the 
maximum principle, to obtain the corresponding discrete extremal flow in the 
states and costates. Treatment in more detail of special cases like the rank 2
case ($n=2$) could also be carried out in the future.  

\section{Conclusions} 

In this paper we have presented the (continuous) geodesic flow on Stiefel 
manifolds with left-invariant metrics. The geodesic equations were obtained 
from two approaches; a variational approach taking reduced variations on 
the Stiefel manifold, and an optimal control approach using costate variables. 
We have attempted to generalize the symmetric representation of the 
$N$-dimensional rigid body flow given by Bloch \etal (2002) to geodesic flows on 
Stiefel manifolds. \\

The solution manifold in the symmetric representation of the rigid body was 
$SO(N)\times SO(N)$. We found it easier to restrict the geodesic flows on Stiefel 
manifolds with a left-invariant metric to the $W^k_{n,N}$, which are symplectic 
submanifolds of $\bbr^{nN}\times\bbr^{nN}$. 
The relation between $SO(N)\times SO(N)$ and $W^k_{N,N}$ for any fixed value of $k$ 
is given by:
\[ SO(N)\times SO(N)= \bigcup_{m\in \{m'\ | \ \frac{k+m'}{2}\in SO(N)\} }W^k_m, \] 
where the union is over all such $m\in\mfrak{so}(N)$ that satisfy $\frac{k+m}{2}
\in SO(N)$. We obtained the geodesic flows using the maximum principle of optimal 
control theory and related them to the Hamiltonian flows using the natural 
symplectic structure on the Stiefel manifold. We also related this optimal control 
formulation to the form naturally derived from variational calculus.\\

Note that the extremal solutions of the optimal control problem in $W^0_{n,N}$ 
are a generalization of the McLachlan-Scovel equations (\ref{macsco}) for the 
$N$-dimensional rigid body, where $Q\T P_0$ is skew symmetric. However, these 
extremal solutions do not generalize the symmetric representation of the rigid 
body equations given in Bloch \etal (2002), wherein $Q\T P$ is orthogonal. 
%\input bangfont %{\bngxii k}

%\subsection*{Acknowledgements}
\ack
The authors would like to thank Yuri Federov and Bo\v{z}idar Jovanovi\'{c} for helpful 
comments and discussions. Support from the NSF is gratefully acknowledged. 

%\begin{thebibliography}{9}
\References

\item[] Bloch, A. M., Crouch, P. E., Marsden, J. E., and Ratiu, T. S., 2002, 
``The symmetric representation of the rigid body equations and their discretization," 
{\em Nonlinearity}, {\bf 15}, 1309-1341. 

\item[] Bloch, A. M., Baillieul, J., Crouch, P. E., and Marsden, J. E., 2003, 
%``Mathematical Preliminaries." In 
{\em Nonholonomic Mechanics and Control}, Vol. 24 of Series in Interdisciplinary 
Applied Mathematics, Springer Verlag, New York. %, 88-92. 

\item[] Bolsinov, A. V., and Jovanovic, B., 2004, ``Complete involutive 
algebras of functions on cotangent bundles of homogeneous spaces," {\em Mathematische 
Zeitschrift}, {\bf 246}, 213-236. 

\item[] Bridges, T., and Reich. S., 2001, ``Computing Lyapunov exponents on a 
Stiefel Manifold," {\em Physica D}, {\bf 156}, 219-238.

\item[] Edelman, A., Arias, T., and Smith, S. T., 1998, ``The Geometry of Algorithms 
with Orthogonality Constraints," {\em SIAM J. Matrix Anal. Appl.}, {\bf 20}, 303-353.

\item[] Elden, L., and Park, H., 1999, ``A Procrustes Problem on the Stiefel 
Manifold," {\em Numerische Mathematik}, {\bf 82}, 599-619.

\item[] Federov, Y. N., 1995, ``Various aspects of $n$-dimensional rigid body 
dynamics," {\em American Mathematical Society Translations}, {\bf 168}, 141-171.

\item[] Federov, Y. N., 2005, ``Integrable flows and Backlund transformations on 
extended Stiefel varieties with application to the Euler top on the Lie group $SO(3)$," 
preprint available at http://arxiv.org/abs/nlin.SI/0505045.

\item[] Gelfand, I. M., and Fomin, S. V., 2000, {\em Calculus of Variations}, 
translated by R. A. Silverman, Dover Publications, Mineola, NY.  

\item[] Kirk, D. E., 2004, {\em Optimal Control Theory: An Introduction}, Dover 
Publications, New York.

\item[] Knorrer, H., 1980, ``Geodesics on the ellipsoid," {\em Invent. Math.}, 
{\bf 59}, 119-143.

\item[] Knorrer, H., 1982, ``Geodesics on quadrics and a mechanical problem of 
C. Neumann," {\em J. Reine Angew. Math.}, {\bf 334}, 69-78.  

\item[] Manakov, S. V., 1976, ``Note on the integration of Euler's equations of the 
dynamics of an $n$-dimensional rigid body," {\em Functional Analysis and its Applications}, 
{\bf 10}, 253-299.

\item[] Marsden, J. E., and Ratiu, T. S., 1999, \emph{Introduction to Mechanics
and Symmetry}, 2nd. ed., Springer-Verlag Inc., New York, 345-348. 

\item[] McLachlan, R. I., and Scovel, C., 1995, ``Equivariant constrained 
symplectic integration," {\em Journal of Nonlinear Science}, {\bf 5}, 233-256. 

\item[] Mischenko, A. S., and Fomenko, A. T., 1978, ``Generalized Liouville 
method of integration of Hamiltonian systems," {\em Functional Analysis and its 
Applications}, {\bf 12}, 113-121. 

\item[] Moser, J., 1980, ``Geometry of quadrics and spectral theory," {\em Chern 
Symposium 1979}, Springer-Verlag Inc., New York, 147-188.

\item[] Moser, J., and Veselov, A. P., 1991, ``Discrete Versions of Some 
Classical Integrable Systems and Factorization of Matrix Polynomials," {\em Communications 
in Mathematical Physics}, {\bf 139}, 217-243. 

\item[] Ratiu, T., 1980, ``The Motion of the Free n-dimensional Rigid Body," 
{\em Indiana University Mathematics Journal}, {\bf 29}, 609-629. 

%\end{thebibliography}
\endrefs

\end{document}